\documentclass[a4paper,11pt]{article}
\usepackage{float} 
\usepackage{fancyhdr}
\usepackage{latexsym,amssymb,amsmath}
\usepackage{amsmath}
\usepackage[french]{babel}
\usepackage[latin1]{inputenc}

\input xy
\xyoption{all}

\newtheorem{defi}{D\'efinition}[section]
\newtheorem{thm}[defi]{Th\'eor\`eme}
\newtheorem{prop}[defi]{Proposition}
\newtheorem{coro}[defi]{Corollaire}
\newtheorem{lemme}[defi]{Lemme}
\newtheorem{rema}[defi]{Remarque}

\newtheorem{nota}[defi]{Notation}
\newtheorem{notas}[defi]{Notations}

\newtheorem{fait}[defi]{Fait}

\newenvironment{dem}{\noindent {\it D\'emonstration $-$ }
                 \noindent}{\hfill $\Box$\vskip 5mm}

\newcommand{\N}{\mathbb{N}}
\newcommand{\Z}{\mathbb{Z}}
\newcommand{\Q}{\mathbb{Q}}
\newcommand{\R}{\mathbb{R}}
\newcommand{\C}{\mathbb{C}}
\newcommand{\K}{\mathbb{K}}

\newcommand{\F}{\mathcal{F}}
\newcommand{\Log}{\mathcal{L}og}
\newcommand{\Pol}{\mathcal{P}ol}
\newcommand{\PP}{\mathcal{P}}

\newcommand{\hh}{\mathcal{H}}
\newcommand{\imp}{\Rightarrow}

\newcommand{\Sym}{\mbox{Sym}}
\newcommand{\Ker}{\mbox{Ker}}
\newcommand{\Hom}{\mbox{Hom}}
\newcommand{\Ext}{\mbox{Ext}}

\newcommand{\isom}{\overset{\sim}{\to}}
 
\newcommand{\OO}{\mathcal{O}}  
\newcommand{\cinf}{\mathcal{C}^{\infty}} 
\newcommand{\aaa}{\mathcal{A}}
\newcommand{\For}{\mbox{For}}

\newcommand{\A}{\mathbb{A}}
\newcommand{\V}{\mathbb{V}}
\newcommand{\G}{\mathcal{G}}

\newcommand{\Eis}{\mathcal{E}is}

\begin{document}

\hfill{2 Mai 2008}
\begin{center}
\vspace{1cm}
{\huge  {\bf R\'ealisation de Hodge du polylogarithme}}
 \\ 
 \vspace{0.3cm}
 {\huge {\bf d'un sch\'ema ab\'elien}} \\
\vspace{1cm}
{\Large David Blotti\`ere}\\
$\;$\\
Universit\"at Paderborn,\\
Institut f\"ur Mathematik, \\
Warburger Str. 100,\\
33098 Paderborn, Allemagne.\\
email: \verb+blottier@math.upb.de+ \\
\vspace{1cm}

\noindent {\Large {\bf R\'esum\'e}} \\

\begin{tabular}{p{15cm}}
 \\
Le r\'esultat principal de cet article est que les courants d\'efinis par Levin dans \cite{l}
permettent de d\'ecrire le polylogarithme d'un sch\'ema ab\'elien au niveau topologique. Ce r\'esultat
avait \'et\'e conjectur\'e par Levin. On en d\'eduit une m\'ethode pour d\'eterminer explicitement
les classes d'Eisenstein des sch\'emas ab\'eliens au niveau topologique. Ces classes ont un int\'er\^et particulier car, d'apr\`es
Kings (cf. \cite{k}), elles ont une origine motivique. Dans \cite{b}, on utilise le r\'esultat principal de cet article
(Corollaire \ref{explicitation_Pol}) pour d\'emontrer que les classes d'Eisenstein des sch\'emas ab\'eliens d'Hilbert-Blumenthal d\'eg\'en\`erent au bord de la compactification de Baily-Borel de la base en une valeur
sp\'eciale de fonction $L$ associ\'ee au corps de nombres totalement r\'eel sous-jacent. On en d\'eduit, dans ce cadre
g\'eom\'etrique, un r\'esultat de non annulation pour certaines classes d'Eisenstein.

\\ \\ \\
\end{tabular}

\noindent {\Large {\bf Abstract}} \\

\begin{tabular}{p{15cm}}
 \\
The main result of this article is the fact that the currents defined by Levin in \cite{l} give
a description of the polylogarithm of an abelian scheme at the topological level. This result
had been conjectured by Levin. This provides a method to explicit the Eisenstein classes of an abelian
scheme at the topological level. These classes are of special interest since they have a motivic origin
by a theorem of Kings (\cite{k}). In \cite{b}, we use the main result of this article
(Corollaire \ref{explicitation_Pol}) to prove that the Eisenstein classes of the universal abelian scheme
over an Hilbert-Blumenthal variety degenerate at the boundary of the Baily-Borel compactification of the base
in a special value of an $L$-function associated to the underlying totally real number field. As a corollary, we get
a non vanishing result for some of these Eisenstein classes in this geometric situation.

\end{tabular}
\end{center}

\newpage

\section{Introduction}

Le polylogarithme de $\mathbb{P}^1 \setminus \{ 0, 1 , \infty\}$ peut \^etre
d\'ecrit explicitement par une matrice dans laquelle apparaissent les logarithmes sup\'erieurs (les fonctions
$Li_k$). 
On peut, gr\^ace \`a cette description, d\'emontrer que les classes d'Eisenstein 
(construites \`a partir du polylogarithme et d'une racine de l'unit\'e)
sont li\'ees aux valeurs sp\'eciales de la fonction $\zeta$ de Riemann. De plus,
elles sont d'origine motivique et engendrent l'image du r\'egulateur.\\

Beilinson et Levin ont d\'efini et d\'ecrit le polylogarithme d'une famille de courbes elliptiques \cite{bl}. 
Pour des courbes elliptiques CM obtenues en tirant par point CM la famille de courbes
elliptiques universelle, les classes d'Eisenstein  (construites \`a partir du polylogarithme et d'une section de torsion)
 sont d'origine motivique et fournissent
un syst\`eme de g\'en\'erateurs de l'image du r\'egulateur (cf. \cite[V-4]{w}).\\

Pour une famille de vari\'et\'es ab\'eliennes, la d\'efinition du polylogarithme
se d\'eduit directement de la th\`ese de Wildeshaus \cite{w}. 
Par analogie avec les deux situations g\'eom\'etriques pr\'ec\'edentes, \'etant donn\'e un sch\'ema ab\'elien de dimension
relative sup\'erieure \`a  $2$, on consid\`ere les questions suivantes: \\
\begin{center}
\begin{tabular}{cp{12cm}}
$(Q_1)$ &  Peut-on d\'ecrire explicitement le polylogarithme ? \\
$(Q_2)$ &  Les classes d'Eisenstein (construites \`a partir du polylogarithme et d'une section de torsion) sont-elles d'origine motivique ? \\
$(Q_3)$ &  Les classes d'Eisenstein sont-elles li\'ees \`a des valeurs sp\'eciales de fonctions $L$  ? \\
$(Q_4)$ & Les classes d'Eisenstein engendrent-elles l'image du r\'egulateur ? \\
\\
\end{tabular}
\end{center}

Dans $\cite{k}$, Kings d\'emontre l'origine motivique des classes d'Eisenstein d'un sch\'ema ab\'elien.\\

Dans \cite{l}, Levin associe \`a un sch\'ema ab\'elien polaris\'e
des courants (nomm\'es courants polylogarithmiques).
Le r\'esultat principal de cet article (Corollaire \ref{explicitation_Pol}) est 
que ces derniers permettent de d\'ecrire le polylogarithmique d'un sch\'ema ab\'elien (au niveau topologique).
Ceci  avait \'et\'e conjectur\'e par Levin. On r\'epond ainsi par l'affirmative \`a la question $(Q_1)$.\\

Dans \cite{b}, on sp\'ecialise la situation aux sch\'emas ab\'eliens d'Hilbert-Blumenthal et 
on utilise ce r\'esultat de fa\c{c}on essentielle pour d\'emontrer que 
les classes d'Eisenstein d\'eg\'en\`erent au bord 
de la compactification de Baily-Borel de la base en une valeur
sp\'eciale de fonction $L$ associ\'ee au corps de nombres totalement r\'eel sous-jacent et en d\'eduire que certaines sont non nulles.
Ainsi, dans cette situation g\'eom\'etrique particuli\`ere, on r\'epond \`a la question $(Q_3)$ par l'affirmative
 et on fait un premier pas dans l'\'etude de la question  $(Q_4)$.\\

On pr\'esente maintenant le contenu de cet article.\\

Dans la section  2, on a rassembl\'e quelques d\'efinitions
et propri\'et\'es concernant les courants. On introduit notamment le complexe des courants 
\`a valeurs dans un fibr\'e vectoriel plat, objet qui intervient dans la formulation de notre r\'esultat principal.\\
Dans la partie suivante, on donne deux d\'efinitions du logarithme d'un sch\'ema ab\'elien;
l'une issue du travail de Wildeshaus (cf. \cite{w}), bas\'ee sur le th\'eor\`eme de Hain-Zucker,
l'autre due \`a Kings (cf. \cite{k}), et on les compare.
On d\'ecrit ensuite le pro-syst\`eme local sous-jacent au logarithme pr\'ec\'edemment d\'efini \`a l'aide du
pro-fibr\'e vectoriel plat construit par Levin (cf. \cite[Part 2]{l}) et on \'enonce les propri\'et\'es du logarithme,
e.g. le r\'esultat du calcul de ses images directes sup\'erieures.\\
Dans la partie 4, on rappelle la d\'efinition du  polylogarithme d'un sch\'ema ab\'elien. Celle-ci fait intervenir
de fa\c{c}on essentielle un morphisme r\'esidu.
Le polylogarithme est une extension de modules de Hodge mixtes qui est rigide, i.e. qui est caract\'eris\'ee par l'extension sous-jacente au niveau topologique.
On d\'emontre qu'il suffit de r\'esoudre une certaine \'equation diff\'erentielle pour expliciter cette derni\`ere
extension (Th\'eor\`eme \ref{pol_equa_diff}). Enfin, les courants de Levin satisfaisant  cette
\'equation diff\'erentielle, on en d\'eduit le r\'esultat principal de ce travail
(Corollaire \ref{explicitation_Pol}). \\
Dans la 
 derni\`ere section, on explique comment on  peut en 
d\'eduire une m\'ethode pour expliciter, au niveau
topologique, les classes d'Eisenstein d'un sch\'ema ab\'elien dont on rappelle auparavant la d\'efinition.

\subsection*{Remerciements}

Ce travail est issu de ma th\`ese de doctorat dirig\'ee par J\"org Wildeshaus. Je tiens \`a le remercier
 pour m'avoir propos\'e ce sujet ainsi que pour les discussions que nous avons partag\'ees.\\
 
 La preuve pr\'esent\'ee ici du Th\'eor\`eme \ref{pol_equa_diff} diff\`ere de celle donn\'ee dans ma th\`ese. Les \'echanges que j'ai eus
 avec Vincent Maillot et Jose I. Burgos \`a propos de la notion de
 courant dans le cadre alg\'ebrique, ainsi que les remarques du rapporteur m'ont permis d'en simplifier
 la d\'emonstration. Je les remercie tous trois. \\
 
  Je remercie également Andrey Levin 
  qui, d'une part a remarqué que l'argument invoqué dans ma thèse pour justifier la lissité des courants polylogarithmiques (cf \cite{l}) 
  était erroné, et d'autre part a eu la gentillesse d'écrire  une preuve de ce résultat dans l'appendice de cet article (cf Proposition A2.1).

\subsection*{Notations et convention}

  Soient $X$ un sch\'ema de type fini, s\'epar\'e et lisse sur $\C$,
$f \colon Y \to Z$ un morphisme entre sch\'emas de type fini, s\'epar\'es et lisses sur $\C$
et $\K \in \{\Q,\C\}$. 
  On note \\

\begin{tabular}{lp{12cm}}
$\overline{X}$ & l'ensemble $X(\C)$ muni de la topologie transcendante, \\
$\overline{f}$ & l'application continue de $\overline{Y}$ vers $\overline{Z}$ induite par $f$,\\
$X^{\infty}$ & la vari\'et\'e diff\'erentielle $\cinf$-r\'eelle associ\'ee \`a $X$, \\
$f^{\infty}$ & l'application lisse de $Y^{\infty}$ vers $Z^{\infty}$ induite par $f$, \\
$\mathcal{F}_{\K}(X)$ & la cat\'egorie des faisceaux en $\K$-vectoriels sur $\overline{X}$, \\
$D_c^b(X)$ & la sous-cat\'egorie pleine de $D^b \mathcal{F}_{\Q}(X)$ ayant pour objets les complexes dont la cohomologie est alg\'ebriquement constructible,\\
$SHM$ & 	la cat\'egorie des $\Q$-structures de Hodge mixtes admissibles polarisables, \\	     
$VSHM(X)\;$   & la cat\'egorie des $\Q$-variations de structures de Hodge mixtes admissibles 
                (cf. \cite{ka}) polarisables
                 sur $X$, \\
$\overline{\mathbb{V}}$ & le (pro-)syst\`eme local sous-jacent \`a $\mathbb{V}$ pour $\mathbb{V} \in Ob((pro\text{-})VSHM(X))$,\\    
$MHM(X)$ & la cat\'egorie des $\Q$-modules de Hodge alg\'ebriques mixtes sur $X$ (cf. \cite{s}).\\	    
         & \\
\end{tabular}

  Par construction, $MHM(X)$ est muni d'un foncteur $rat$ de $MHM(X)$ vers $Perv(X)$, 
le coeur de la $t$-structure perverse autoduale sur $D_c^b(X)$, qui est fid\`ele et exact.
  Celui-ci induit  un foncteur de $D^bMHM(X)$  vers $D^bPerv(X)$ qui compos\'e avec le 
foncteur $real$ de Beilinson (cf. \cite{bbd}) fournit un foncteur d'oubli 
$\For \colon D^b MHM(X) \to  D_c^b(X)$. On dispose \'egalement d'un foncteur
$\iota_X \colon VSHM(X) \to MHM(X)$ qui est exact, pleinement fid\`ele et gr\^ace auquel on identifie
$VSHM(X)$ \`a une sous-cat\'egorie pleine de $MHM(X)$. 
  Le foncteur $\For$ associe \`a un objet de $VSHM(X)$ le syst\`eme local sous-jacent d\'ecal\'e.
  Dans ce texte, on fait la convention suivante:
{\it l'image d'un objet de $VSHM(X)$ sous $For$ est son syst\`eme local sous-jacent concentr\'e en degr\'e $0$}, 
i.e. on ne tient pas compte du d\'ecalage. \\

 On fixe $i$ une racine carr\'ee de $-1$ dans $\C$ pour la suite.
Ce choix d\'etermine une orientation canonique des vari\'et\'es diff\'erentielles r\'eelles
associ\'ees aux sch\'emas de type fini, s\'epar\'es et lisses sur $\C$.\\

Soit $S$ un sch\'ema de type fini, s\'epar\'e, connexe et lisse sur $\C$ et soit:\\

\begin{tabular}{cp{14cm}}

$\pi \colon A \to S $ &  un sch\'ema ab\'elien de section unit\'e $e$ et de dimension relative pure $d$, \\
$\hh$ & $:=(R^1 \pi_* \Q)^{\vee}$, objet pur de poids $-1$ de $VSHM(S)$, \\
$j \colon U \hookrightarrow A$ &  l'immersion ouverte compl\'ementaire de $e$, \\
          $\pi_U$ & $:= \pi \circ j$.

\end{tabular}

\section{Courants}

\subsection{Courants sur une vari\'et\'e diff\'erentielle}

Soit $X$ une variété différentielle de dimension pure $n$.
\begin{notas} $-$  Soit $p \in \N$, $0 \leq p \leq n$. On note: \\

\begin{tabular}{cl}
 $\OO_X$ & le faisceau des fonctions diff\'erentielles sur $X$ à valeurs dans
          $\C$,\\
 $\Omega_X^p$ & le faisceau des $p$-formes diff\'erentielles complexes de $X$, \\
 $\Omega_{X,c}^p$ &  le faisceau des $p$-formes diff\'erentielles complexes
                        \`a supports compacts de $X$.
\end{tabular}
\end{notas}

On munit  $\Omega_X^p(X)$ de la topologie donn\'ee par [D, 17.2].
Pour $K \subset X$ compact,
l'espace des $p$-formes diff\'erentielles complexes sur $X$ \`e support dans $K$, not\'e $\Omega^p_X (X,K)$,
h\'erite de la topologie induite, qui en fait un espace de Fr\'echet.

\begin{defi} $-$
Un $p$-courant sur $X$ est une forme lin\'eaire $$T : \Omega^{n-p}_{X,c}(X) \to \C$$ dont
la restriction \`a chacun des $\Omega_X^p(X,K)$ ($K \subset X$ compact) est continue.
On note $\aaa_X^{p}(X)$ l'espace des $p$-courants sur $X$.
\end{defi}

On munit $\aaa_X^{p}(X)$ de la topologie faible qui est induite 
par les semi-normes
   $$ T \in \aaa_X^{p}(X) \mapsto | T(\alpha) | $$
pour $\alpha \in \Omega^{n-p}_{X,c}(X)$ (cf. [D, 17.8]).\\

Soient $U,V$ deux ouverts de $X$, $U \subset V$, et $K \subset U$ un compact de $K$.
On a une application naturelle
 $\Omega_X^{n-p}(U,K) \to \Omega_X^{n-p}(V,K)$ (prolongement par $0$ sur $V \setminus U$). On en d\'eduit
une application de restriction
$$res^V_U :  \aaa^{p}_X(V):= \aaa^{p}_V(V) \to \aaa^{p}_U(U) =: \aaa^{p}_X(U).$$
On d\'efinit ainsi un pr\'efaisceau sur $X$ not\'e $\aaa^{p}_X$. Ce pr\'efaisceau est un faisceau
(cf. [D, 17.4.2]).\\

Soient $p,q \in \N$ tels que $p+q \leq n$. On a un accouplement canonique:
$$ \begin{array}{rccccc}
   \psi_{p,q} :  & \aaa^{p}_X \otimes  \Omega^q_X & \to & \aaa^{p+q}_X .&  \\
                 &   T \otimes \omega & \mapsto & T(\omega \wedge \cdot ) & \end{array}$$
En particulier, $\psi_{p,0}$ d\'efinit une structure de $\OO_X$-module sur $\aaa^{p}_X$. \\

\subsection{Courants sur une vari\'et\'e diff\'erentielle orient\'ee}

Supposons que $X$ est orient\'ee.
Soient $U \subset X$ ouvert et $p \in \N$,
$0 \leq p \leq n$.
On dispose alors de l'int\'egrale
$$ \int_U : \Omega^n_{X,c}(U) \to \C $$
gr\^ace \`a laquelle, \`a $\eta \in \Omega^p_X(U)$, on associe un $p$-courant sur $U$ not\'e $T_{\eta}$ d\'efini
par:
$$ \begin{array}{rccccc}
    T_{\eta} & : &  \Omega^{n-p}_{X,c}(U) & \to & \C & .\\

               &   &  \omega & \mapsto & \int_U \eta \wedge \omega & \end{array} $$
L'association $\eta \in \Omega^p_X(U) \mapsto T_{\eta}  \in \aaa^{p}_X(U) $ donne un monomorphisme de faisceaux
not\'e $Int_p$. Dans la suite, on notera simplement  $\eta$ le courant $T_{\eta}$.
Les $p$-courants sur $X$ qui viennent d'une $p$-forme diff\'erentielle sur $X$,
via $Int_p(X)$, sont appel\'es courants lisses. \\

\subsection{Courant associ\'e \`a une sous-vari\'et\'e ferm\'ee orient\'ee}

Soit $i : Y \hookrightarrow X$ une immersion ferm\'ee. On suppose que $Y$ est orient\'ee et
on note $m$ sa dimension suppos\'ee pure. Alors, l'application
$$
\begin{array}{ccl}
\Omega_{X,c}^m(X) & \to & \C \\
     \omega & \mapsto & \int_Y i^* \omega
\end{array}
$$
d\'efinit un $(n-m)$-courant que l'on note $\delta_Y$.

\subsection{Diff\'erentiation des courants} \label{derivcourants}

Soit $T$ un $p$-courant sur $X$. On d\'efinit la diff\'erentielle de $T$, 
not\'ee $dT$, comme \'etant le $(p+1)$-courant d\'efini par:
$$ dT(\omega) = (-1)^{p+1}  T (d \omega) $$
pour $\omega \in \Omega^{n-p-1}_{X,c}(X)$. 
Le facteur $(-1)^{p+1}$ est ajout\'e pour que la diff\'erentiation des courants soit compatible avec celle des formes diff\'erentielles. En effet, avec la pr\'ec\'edente d\'efinition, si $X$ est orient\'ee, on a:
$$ d \circ Int_p(\eta) = Int_{p+1} \circ d(\eta) $$
pour $\eta \in \Omega_X^p(X)$.\\

\begin{lemme} $-$ \label{lincourants}
Soient $p,q \in \N$, $0 \leq p \leq n-1$, $0 \leq q \leq p + 1$. 
Soient $T$ un $p$-courant sur $X$ et 
$\omega \in \Omega^q_X(X)$. Alors on a: \\

\begin{itemize}
\item[1.]  $d \circ d (T) = 0.$ 
\item[2.]  $d \psi_{p,q} (T \otimes \omega)  
            = \psi_{p+1,q} (dT \otimes \omega) + 
              (-1)^p \; \psi_{p,q+1}(T \otimes d \omega)$.
\end{itemize}
\end{lemme}

\begin{dem} La premi\`ere \'egalit\'e se d\'eduit de la propri\'et\'e $d \circ d (\eta)= 0$ pour tout $\eta \in \Omega^{n-p-1}_{X,c}(X)$. Pour prouver
 la deuxi\`eme, on consid\`ere $\eta \in  \Omega^{n-p-q-1}_{X,c}(X)$ et on effectue le calcul suivant pour conclure.

$$ \begin{array}{lll} 
   d \psi_{p,q} ( T \otimes \omega ) (\eta) & = & 
                (-1)^{p+q+1} \; T ( \omega \wedge d\eta) \\
                 & = & (-1)^{p+1}  \; 
                 T ( d(\omega \wedge \eta) - d \omega \wedge \eta ) \\
                 & = & dT (\omega \wedge \eta) + 
                       (-1)^{p+2}  \; T( d\omega \wedge \eta ) \\
                 & = &  \psi_{p+1,q}(dT \otimes \omega) (\eta) 
                    + (-1)^{p} \; \psi_{p,q+1}(T \otimes d\omega) (\eta).
 \end{array}
$$

\end{dem}

\subsection{Courants \`a valeurs dans un fibr\'e vectoriel}

Soient $E$ un fibr\'e vectoriel
complexe de rang $N$ au-dessus de $X$ et $p \in \N$, $0 \leq p \leq n$.

\begin{defi} $-$
Le faisceau des $p$-courants sur $X$ \`a valeurs dans $E$ est 
    $$\aaa^p_{X}(E) := \aaa^p_{X} \otimes_{\OO_X} E.$$
\end{defi}

Comme $\aaa^p_{X}(E)$ est un $\OO_X$-module, le faisceau $\aaa^p_{X}(E)$ est fin.
Dans le cas o\`u la vari\'et\'e est orient\'ee, on a un monomorphisme de faisceaux
$$ Int_p \otimes \text{Id}_E \colon \Omega^p_X(E) \to \aaa^p_{X}(E) . $$
Un  $p$-courant sur $X$ à valeurs dans $E$ est dit lisse s'il provient, via 
$ Int_p \otimes \text{Id}_E$, d'une  $p$-forme diff\'erentielle sur $X$ \`a valeurs dans $E$.

\subsection{Notion de convergence} \label{convergencecourants}
On cherche \`a d\'efinir une notion de convergence pour les courants
\`a valeurs dans un fibr\'e vectoriel.

\subsubsection{Cas o\`u le fibr\'e est trivial}
Si $E$ est le fibr\'e trivial de rang $N$ sur $X$, alors on a
la d\'ecomposition
$$ \aaa^p_{X}(E) =  (\aaa^p_{X})^N $$
relativement \`a la base canonique de $\C^N$ not\'ee $(e_1,..,e_N)$ 
et on a une notion naturelle de convergence sur $\Gamma(X,\aaa^p_{X}(E))$.
En effet,
soit $(T_k)_{k \geq 0}$ une suite d'\'el\'ements de $\Gamma(X,\aaa^p_{X}(E))$ et 
$T \in  \Gamma(X,\aaa^p_{X}(E))$.
Pour tout $k \geq 0$, on écrit
 $$ T_k = \sum_{1 \leq i \leq N} \; T_k^i e_i $$
la d\'ecomposition de $T_k$ relativement \`a la base canonique de $\C^N$.
On d\'ecompose de m\^eme $T$,
 $$ T= \sum_{1 \leq i \leq N} \; T^i e_i .$$

\begin{defi}{} $-$ \label{convtriv} Dans cette situation,
on dit que $(T_k)_{k \geq 0}$ tend vers $T$ dans $\Gamma(X,\aaa^p_{X}(E))$
et on écrit
 $ T_k \underset{k \to \infty}{\to} T $
si pour tout
$( \omega_1 ,.. , \omega_N)  \in (\Omega_{X,c}^{n-p}(X))^N$
$$ \left( T^1_k(\omega_1) , ..,  T^N_k(\omega_N) \right)  \underset{k \to \infty}{\to}
   ( T^1(\omega_1) , ..,  T^N(\omega_N)) \mbox{ dans } \C^N. $$
\end{defi}

La notion de convergence de la D\'efinition \ref{convtriv} est invariante par automorphisme, comme on le v\'erifie ci-dessous.
Soit $\varphi : E \to E$
un automorphisme de fibr\'e vectoriel donn\'e relativement \`a la base canonique de $\C^N$, par
$$
\begin{array}{ccl}
X & \to & GL_N(\C) \\
x & \mapsto & (\varphi_{ij}(x))_{1 \leq i , j \leq N}
\end{array}
$$
où $\varphi_{ij} \in \OO_X(X)$.
Alors, $\varphi$ induit un isomorphisme
$$ Id \otimes \varphi^* : \Gamma(X,\aaa^p_{X}(E)) \isom \Gamma(X,\aaa^p_{X}(E)) $$
qu'on explicite.
Si $T = \displaystyle \sum_{1 \leq i \leq N} \; T^i e_i $, alors
$$
\begin{array}{ccl}
  Id \otimes \varphi^* (T) & = &  \displaystyle \sum_{1 \leq i \leq N} \sum_{1 \leq j \leq N}
                T^i \otimes \varphi_{ji} e_j \\
                             & = &
        \displaystyle  \sum_{1 \leq i \leq N} \sum_{1 \leq j \leq N}   \varphi_{ji} T^i \otimes  e_j \\
 & = &
        \displaystyle   \sum_{1 \leq j \leq N}    \left( \sum_{1 \leq i \leq N}\varphi_{ji} T^i \right ) \otimes  e_j .
\end{array}
$$
De cette formule, on d\'eduit le

\begin{lemme} $-$ \label{colleconv} Etant donn\'es $(T_k)_{k \geq 0}$ une
suite d'\'el\'ements de $\Gamma(X,\aaa^p_{X}(E))$ et
$T \in  \Gamma(X,\aaa^p_{X}(E))$, on a l'\'equivalence:
$$ T_k \underset{k \to \infty}{\to} T  \Longleftrightarrow
   Id \otimes \varphi^* (T_k) \underset{k \to \infty}{\to} Id \otimes \varphi^* (T) .$$
\end{lemme}

La notion de convergence de la D\'efinition \ref{convtriv} est locale. En effet, on a le 

\begin{lemme} $-$ \label{locconv} Soient  $(T_k)_{k \geq 0}$ une
suite d'éléments de $\Gamma(X,\aaa^p_{X}(E))$ et $T \in  \Gamma(X,\aaa^p_{X}(E))$, soit 
$(U_i)_{i \in I}$ un recouvrement ouvert de $X$.  On a:

$$ T_k \underset{k \to \infty}{\to} T  \Longleftrightarrow \left(  \forall i \in I \quad
   res^X_{U_i}(T_k) \underset{k \to \infty}{\to} res^X_{U_i}( T)  \right).$$
\end{lemme} 

\begin{dem} L'implication $\imp$ est triviale. Pour d\'emontrer l'autre, il suffit d'utiliser une partition de l'unit\'e
adapt\'ee au recouvrement $(U_i)_{i \in I}$.
\end{dem}

\subsubsection{Cas g\'en\'eral}

\begin{defi}  $-$  Soit $E$ un fibr\'e vectoriel complexe de rang $N$ sur $X$. Une famille 
$(U_i,\varphi_i)_{i \in I}$ o\`u \\

\begin{itemize}
\item[a)] $(U_i)_{i \in I}$ est un recouvrement ouvert de $X$, 
\item[b)] pour tout $i \in I$, $\varphi_i$ est
          un isomorphisme de fibr\'es vectoriels,
$$
\xymatrix{ U_i \times \C^N \ar[r]^-{\varphi} \ar[d]^-{pr_1} & E_{U_i} \ar[dl] \\
           U_i }$$
\item[]
\end{itemize}est appel\'ee famille de trivialisations locales de $E$.
\end{defi}

On \'etend la D\'efinition \ref{convtriv} comme suit.

\begin{defi}  \label{defconv} $-$ Si $E$ est un fibr\'e vectoriel complexe de rang $N$ sur $X$, si
$(T_k)_{k \geq 0}$ est une
suite d'\'el\'ements de $\Gamma(X,\aaa^p_{X}(E))$ et
$T \in  \Gamma(X,\aaa^p_{X}(E))$, on dit que
$$ (T_k)_{k \geq 0} \mbox { tend vers } T \mbox{ dans }  \Gamma(X,\aaa^p_{X}(E)) $$
si et seulement si, il existe une famille de trivialisations locales de $E$  $(U_i,\varphi_i)_{i \in I}$ telle que
pour tout $i \in I$,

$$  Id \otimes \varphi_i^* (res^X_{U_i}(T_{k}))   \underset{k \to \infty}{\to}   Id \otimes \varphi_i^* (res^X_{U_i}(T))
   \mbox{ dans } \Gamma(U_i,\aaa^p_{U_i}(U_i \times \C^N)) $$
au sens de la D\'efinition \ref{convtriv}.
\end{defi}

\begin{rema} $-$ Les Lemmes \ref{colleconv} et \ref{locconv} assurent que cette d\'efinition est compatible
avec la d\'efinition \ref{convtriv} et que, lorsque la condition de convergence vaut pour une famille de trivialisations
 locales, elle vaut pour toutes.
\end{rema}

\subsection{Complexe des courants associ\'e \`a un fibr\'e vectoriel plat} \label{courantsplat}
Soit $E$ un fibr\'e vectoriel
complexe de rang $N$ au-dessus de $X$ muni d'une connexion plate
$\nabla :E \to \Omega^1_{X} \otimes E$. On note
$(\Omega^{\bullet}_{X} \otimes E ,\nabla^{\bullet})$ le complexe de de Rham de $(E,\nabla)$.

Soit $\nabla'^p \colon \mathcal{A}_X^p \otimes E 
      \to \mathcal{A}_X^{p+1}\otimes E$ l'unique morphisme de faisceaux 
de $\mathcal{O}_X$-modules caract\'eris\'e par la condition suivante. Pour tout $U$ ouvert simplement connexe de $X$, 
$T \in \Gamma(U,\mathcal{A}^p_X)$, 
$s \in \Gamma(U,E)$ avec 
$$\nabla(s) = \sum_{i \in I_s}  \omega_i \otimes  s_i$$
o\`u $\omega_i \in \Gamma(U,\Omega^1_X)$  et $s_i \in  \Gamma(U,E)$ pour tout $i$ dans l'ensemble d'indices $I_s$:
$$ \nabla'^p ( T \otimes s ) = 
    dT \otimes s + (-1)^p \sum_{i \in I_s} \psi_{1,p}(T \otimes \omega_i) \otimes s_i.$$

On peut v\'erifier que ce morphisme est bien $\mathcal{O}_X$-lin\'eaire en chacune des deux composantes au moyen du 2. du lemme \ref{lincourants}.

\begin{lemme}{} \label{bord}
$-$ Pour $p \in \N$, $0 \leq p \leq n-1$, $\nabla'^{p+1} \circ \nabla'^p = 0$.
\end{lemme}
\noindent 

\begin{dem}
L'assertion est de nature locale. On peut donc supposer que $E$ est le fibr\'e trivial de rang $N$ et que $\nabla$ est la connexion de Gau\ss{}-Manin. Il suffit en fait de consid\'erer le cas $N=1$. Mais alors, l'assertion r\'esulte du 1. du lemme \ref{lincourants}.
\end{dem}

\begin{defi}{} $-$
Le complexe des courants sur $X$ \`a valeurs dans $E$ est le complexe
 $$[\dots \to 0 \to \underset{deg. \; 0}{\aaa^{0}_{X}(E)} \overset{\nabla'^0}{\to}
           \underset{deg. \; 1}{\aaa^{1}_{X}(E)} \overset{\nabla'^1}{\to}
           \underset{deg. \; 2}{\aaa^{2}_{X}(E)} \overset{\nabla'^2}{\to} \dots
           \overset{\nabla'^{n-1}}{\to} \underset{deg. \; n}{\aaa^{n}_{X}(E)} \to 0 \to \dots].$$
\end{defi}

On suppose maintenant que $X$ est orient\'ee.

\begin{prop}{} $-$ 
Le morphisme de complexes
$$\xymatrix{ [\dots \ar[r] & 0 \ar[r] \ar@{=}[d] &
                    \Omega^{0}_{X}(E)  \ar[r]^-{\nabla^{0}} \ar[d]^-{Int_0 \otimes Id} &
                    \Omega^{1}_{X}(E)  \ar[r]^-{\nabla^{1}} \ar[d]^-{Int_1 \otimes Id} &
                \dots \ar[r]^-{\nabla^{n-1}}  & \Omega^{n}_{X}(E) \ar[d]^-{Int_{n} \otimes Id } \ar[r] &
                   0 \ar[r] \ar@{=}[d] & \dots ] \\
             [\dots \ar[r] & 0 \ar[r]   & 
                    \aaa^{0}_{X}(E)  \ar[r]^-{\nabla'^{0}}  &
            \aaa^{1}_{X}(E)  \ar[r]^-{\nabla'^{1}} & 
                \dots \ar[r]^-{\nabla'^{n-1}}  & \aaa^{n}_{X}(E)  \ar[r] &
             0 \ar[r]  & \dots ]}$$
est un quasi-isomorphisme.
\end{prop}

\begin{dem}
L'assertion est de nature locale.
Il suffit de prouver le r\'esultat pour $X$ une boule ouverte
de $\R^n$ et $E=X \times \C^N \overset{pr_1}{\to} X$
le fibr\'e trivial au-dessus de $X$ muni de la connexion de Gauss-Manin $\nabla_{GM}$.
On se ram\`ene alors au cas $N=1$. Pour la preuve du r\'esultat dans cette situation, on renvoie
\`a [GH, p. 382]. \\

\end{dem}

\begin{coro}{} \label{resolution} $-$ La suite
$$ \xymatrix{
   0 \ar[r] &  Ker(\nabla) \ar[r]^-{Int_0 \otimes Id} &
     \aaa^{0}_{X}(E)  \ar[r]^-{\nabla'^{0}}  &
            \aaa^{1}_{X}(E)  \ar[r]^-{\nabla'^{1}} &
                \dots \ar[r]^-{\nabla'^{n-1}}  & \aaa^{n}_{X}(E)  \ar[r] &
             0}$$
est une suite exacte longue.
\end{coro}

\section{Le logarithme d'un sch\'ema ab\'elien}

\subsection{D\'efinition issue de la th\`ese de Wildeshaus \cite{w}}

\subsubsection{Cas absolu} \label{deflogabsolu}

\begin{nota} $-$ On d\'esigne par $VSHMU(A)$ la sous-cat\'egorie pleine de $VSHM(A)$ dont les
             objets sont les variations unipotentes, i.e. qui admettent une filtration dont les gradu\'es
             sont des variations constantes.
\end{nota}

On suppose dans cette partie que $S = \text{Spec}(\C)$ et on fixe $a \in A(\C)$.
La $\Q$-alg\`ebre $\Q[\pi_1(\overline{A},a)]$ est munie d'une augmentation canonique
$\varepsilon \colon \Q[\pi_1(\overline{A},a)] \to \Q$ dont on note $\mathfrak{a}_a$ le noyau.\\

La th\'eorie des int\'egrales it\'er\'ees de Chen permet de munir 
chacun des $\Q[\pi_1(\overline{A},a)] / \mathfrak{a}_a^n$ de $\Q$-structures de Hodges mixtes
canoniques pour $n \in \N^*$. De plus, les morphismes de projection
$$ pr_{n,m} \colon \Q[\pi_1(\overline{X},x)] / \mathfrak{a}^n \to \Q[\pi_1(\overline{X},x)] /
   \mathfrak{a}^m, \quad m,n \in \N^*, \; m \leq n$$
sont des morphismes de $\Q$-structures de Hodge. 
On dispose ainsi d'une pro-$\Q$-structure de Hodge mixte sur le pro-$\Q$-vectoriel 
$$\Q[\pi_1(\overline{A},a)]^{\; \widehat{\;}} := 
   \displaystyle \lim_{\overset{\longleftarrow}{n \geq 1}} \; Q[\pi_1(\overline{A},a)]
   / \mathfrak{a}_a^n,$$ 
o\`u les morphismes de transitions de la limite projective sont les projections
$pr_{n,m}$, qui est telle que: \\

\begin{tabular}{cp{14.8cm}} 
a) & le morphisme de structure de $\Q$-algèbre , 
        $\Q \to \Q[\pi_1(\overline{A},a)]^{\; \widehat{\;}}$ 
     est sous-jacent à un morphisme de pro-structures de Hodge mixtes $1 : \Q(0) \to 
     \Q[\pi_1(\overline{A},a)]^{\; \widehat{\;}}$, \\
b) & la multiplication dans  $\Q[\pi_1(\overline{A},a)]^{\; \widehat{\;}}$ est
     un morphisme de la cat\'egorie pro-$SHM$,\\
c) & pour $\mathbb{V} \in Ob(VSHMU(A))$, la repr\'esentation de monodromie
     $\pi_1(\overline{A},a) \to \text{End}(\overline{\mathbb{V}_a})$ induit un morphisme
     de la cat\'egorie pro-$SHM$ $\rho_a \colon  
     \Q[\pi_1(\overline{A},a)]^{\; \widehat{\;}} \to \underline{\text{End}}(\mathbb{V}_a)$. \\
     \\
\end{tabular}        

On peut alors rappeler l'\'enonc\'e du th\'eor\`eme de Hain-Zucker. 

\begin{thm} \emph{\cite[Thm 1.6]{hz}} $-$ 
Le foncteur 
$$\begin{array}{ccc} VSHMU(A) & \to & 
 \left( \begin{array}{c} V \in Ob(SHM) \; \mbox{ muni d'un morphisme de } \mbox{pro-}SHM \\
                         \Q[\pi_1(\overline{A},a)]^{\; \widehat{\;}} \to 
                         \underline{\emph{End}}(V)  
\end{array} \right) \\
& & \\
\mathbb{V} & \mapsto & (\mathbb{V}_a,\rho_a) \end{array} $$
est une équivalence de catégories.
\end{thm}

\begin{defi} $-$ On applique ce théorème
à $\Q[\pi_1(\overline{A},a)]^{\; \widehat{\;}}$ muni de la représentation
donnée par la multiplication. On obtient un objet de pro-$VSHMU(A)$, 
le logarithme de $A$ que l'on note $\Log_{A,a}$ ou simplement $\Log_a$ lorsqu'il n'y a pas d'ambigu\"it\'e
sur la vari\'et\'e ab\'elienne consid\'er\'ee.
\end{defi}

Le logarithme est en outre caract\'eris\'e par la propri\'et\'e universelle suivante.

\begin{thm} $-$ \label{caraclog} Le foncteur  
$$\begin{array}{ccc} VSHMU(A) & \to & Ab \\
\mathbb{V} & \mapsto &\mbox{\emph{Hom}}_{SHM}(\Q(0), \mathbb{V}_a) \end{array}$$
est pro-représenté par $\Log_{a}$, i.e. on a une bijection naturelle:
$$ \mbox{\emph{Hom}}_{\mbox{\scriptsize pro-}VSHMU(A)}(\Log_{a},\mathbb{V}) \to 
   \mbox{\emph{Hom}}_{SHM}(\Q(0), \mathbb{V}_a) \quad , \quad 
   \varphi \mapsto \varphi_a \circ 1.$$
\end{thm}

Cet \'enonc\'e est \'equivalent au th\'eor\`eme de Hain-Zucker.
Le pro-système local sous-jacent à $\Log_{A,a}$ est lui aussi caractérisé par une propriété
universelle.

\begin{thm} $-$  \label{caraclogtop} Le foncteur 
$$\begin{array}{ccccc} 

\F_{\Q}(A)  &  \underset{\text{pleine}}{\supseteq} & \left( \begin{array}{c} \Q \text{-syst\`emes locaux sur } \overline{A} \text{ admettant }\\
                        \text{  une filtration dont les gradu\'es } \\
			\text{  sont des faisceaux constants } \\
	\end{array} \right)
&  \to & \Q\mbox{-}vect \\
& & \mathbb{V} & \mapsto  & \mathbb{V}_a \end{array}$$
est pro-représenté par $\overline{\Log_{a}}$, i.e. on a un bijection naturelle:
$$ \mbox{\emph{Hom}}_{\mbox{\scriptsize pro-}\F_{\Q}(A)}(
    \overline{\Log_{a}},\mathbb{V}) \to 
     \mathbb{V}_a \quad , \quad 
   \varphi \mapsto \varphi_a ( 1).$$
\end{thm}

\subsubsection{Cas relatif} \label{deflogrel}

\begin{nota} $-$ Soit $VSHMU(A,\pi)$ la sous-cat\'egorie pleine de $VSHM(A)$ dont les
             objets sont les variations unipotentes relativement \`a $\pi$, i.e. 
	     qui admettent une filtration dont les gradu\'es sont dans l'image de
	     $ \pi^* \colon VSHM(S) \to VSHM(A)$.
\end{nota}

Soit $s \in S(\C)$ et $a:=e(s) \in \A(\C)$. D'apr\`es le th\'eor\`eme de 
Ehresmann, $\overline{\pi}$ est une fibration
localement triviale et donc on a la suite exacte scind\'ee suivante:
$$ \xymatrix{ 1 \ar[r] & \pi_1(\overline{A}_s,a) \ar[r] & \pi_1(\overline{A},a) \ar[r]_{\overline{\pi}_*} 
 & \pi_1(\overline{S},s)  \ar@/_1pc/[l]_{\overline{e}_*} \ar[r] & 1 } \; .$$
On d\'efinit une action de $\pi_1(\overline{A},a)  = \pi_1(\overline{A}_s,a)
 \rtimes  \pi_1(\overline{S},s)$
sur $\Q[\pi_1(\overline{A}_s,a)] {\;}^{\widehat{\;}}$ en faisant
agir $\pi_1(\overline{A}_s,a)$ par multiplication \`a gauche et $\pi_1(\overline{S},s)$ par conjugaison. 
On a ainsi construit un pro-syst\`eme local de $\Q$-vectoriels sur $\overline{A}$ que l'on note 
$\mathbb{V}$. On v\'erifie que la fibre en $s' \in S(\C)$ de $\mathbb{V}$  
s'identifie canoniquement \`a $\overline{\Log_{A_{s'},e(s')}}$. Ainsi, fibre \`a fibre, $\mathbb{V}$ est muni d'une filtration par le poids et d'une filtration de Hodge d'apr\`es la cas absolu pr\'ec\'edemment trait\'e. 

\begin{thm} \emph{\cite[I-Thm 3.3]{w}} $-$ Le pro-syst\`eme local $\mathbb{V}$ muni de ces deux filtrations d\'efinies fibre \`a fibre d\'efinit 
un objet de pro-$VSHMU(A,\pi)$.
\end{thm}

\begin{defi} $-$ L'objet de pro-$VSHMU(A,\pi)$ du th\'eor\`eme pr\'ec\'edent est appel\'e logarithme de 
$A/S$ et est not\'e $\Log_{A/S,s}$ ou simplement $\Log_{s}$ lorsqu'il n'y a pas de confusion possible quant au sch\'ema ab\'elien que l'on consid\`ere.
\end{defi}

On a un morphisme de variations de $\Q$-structures de Hodge sur $S$ canonique,
$1 : \Q(0) \to e^ * \Log_s$ qui est induit par la structure de $\Q$-alg\`ebre
de $\Q[\pi_1(\overline{A}_s,a)] {\;}^{\widehat{\;}}$. On caract\'erise maintenant
$\Log_s$ et $\overline{\Log_s}$ par les propri\'et\'es suivantes qui sont des versions
relatives des Th\'eor\`emes \ref{caraclog} et \ref{caraclogtop}.

\begin{thm} \emph{\cite[I-Thm 3.5]{w}}   $-$ \label{thmutile}
La transformation naturelle entre foncteurs de  $VSHMU(A,\pi)$ vers $VSHM(S)$:
$$ \begin{array}{ccc}  
   \pi_* \underline{\mbox{Hom}}(\Log_s, \cdot) & \to & e^* \\
           \varphi & \mapsto & e^* (\varphi)(1) \end{array} $$
est un isomorphisme de foncteurs. \\
\end{thm}

\begin{thm} \emph{\cite[I-Thm 3.5]{w}}  $-$
La transformation naturelle entre foncteurs de la cat\'egorie  
 $$\left( \begin{array}{c} \Q \text{-syst\`emes locaux sur } \overline{A} 
                             \text{ admettant  une filtration dont les }\\
                        \text{   gradu\'es sont des pullbacks par } \pi
			  \text{ de syst\`emes} \text{ locaux sur } \overline{S} \\
	\end{array} \right) \underset{\text{pleine}}{\subseteq} \F_{\Q}(A) $$ vers 
la cat\'egorie des $\Q$-syst\`emes locaux sur $\overline{S}$:
$$ \begin{array}{ccc}  
   \overline{\pi}_* \underline{Hom}(\overline{\Log_{s}}, \cdot) & \to & \overline{e}^* \\
           \varphi & \mapsto & (\overline{e}^* \varphi)(\overline{1}) \end{array} $$
est un isomorphisme de foncteurs. \\
\end{thm}

Soit $s' \in S(\C)$. Tout chemin allant de $s'$ à $s$ induit un isomorphisme de 
	pro-variations $\Log_{s'} \isom \Log_{s}$. La propriété universelle du logarithme
	implique que cet isomorphisme est en fait indépendant du choix de chemin.
	Ainsi, on note simplement $\Log$ l'objet $\Log_{s}$. On pourra \'egalement noter
	le logarithme $\Log_{A/S}$ lorsque l'on voudra pr\'eciser le sch\'ema ab\'elien.

\subsection{Comparaison avec la d\'efinition du logarithme due \`a Kings \cite{k}}

Soient $s \in S(\C)$ et $a := e(s)$.
On note $\mathfrak{a}_s$ le noyau de l'augmentation
$\varepsilon_s : \Q[\pi_1(\overline{A}_s,a)] \to \Q$.
On a une suite exacte de $\Q$-vectoriels munie d'un scindage canonique:
$$\xymatrix{  (\Sigma) \quad  \quad \quad  & 0 \ar[r]
& \mathfrak{a}_s / \mathfrak{a}_s^2 \ar[r] &
 \Q[\pi_1(\overline{A}_s,a)]/ \mathfrak{a}_s^2 \ar[r]^-{\overline{\varepsilon^{(1)}_s}} &
\Q \ar@<1ex>[l]^-{\overline{1^{(1)}_s}} \ar[r] & 0}.$$
dans laquelle  $\varepsilon^{(1)}_s$ (resp. $1^{(1)}_s$) est l'augmentation (resp. le morphisme
de structure de $\Q$-alg\`ebre) 
de $\Q[\pi_1(\overline{A}_s,a)]/ \mathfrak{a}_s^2$.
On munit chacun des termes de cette suite d'une action de
$\pi_1(\overline{A},a)= \pi_1(\overline{A}_s,a) \rtimes  \pi_1(\overline{S},s)$.
Le groupe $\pi_1(\overline{A},a)$
agit trivialement sur $\Q$.
Sur $\Q[\pi_1(\overline{A}_s,a)]/ \mathfrak{a}_s^2$ et $\mathfrak{a}_s / \mathfrak{a}_s^2$,
$\pi_1(\overline{A}_s,a)$ agit par multiplication et $\pi_1(\overline{S},s)$ agit par conjugaison.
On remarque que $\pi_1(\overline{A}_s,a)$  agit trivialement sur $\mathfrak{a}_s / \mathfrak{a}_s^2$.
 On vérifie
alors que $(\Sigma)$ est une suite exacte de $\pi_1(\overline{A},a)$-modules (non scind\'ee si $d \geq 1$),
que l'on consid\`ere comme une suite exacte de $\Q$-systèmes locaux sur $\overline{A}$. \\

On va maintenant installer des filtrations sur ces sytèmes locaux. Pour tout $s' \in S(\C)$,
on applique le foncteur \guillemotleft restriction à $\overline{A}_{s'}$\guillemotright\;
 à $(\Sigma)$. Le résultat est une
suite exacte de $\pi_1(\overline{A}_{s'},e(s'))$-modules 
canoniquement isomorphe à:
$$ 0 \to  \mathfrak{a}_{s'} / \mathfrak{a}_{s'}^2  \to
    \Q[\pi_1(\overline{A}_{s'},e(s'))] / \mathfrak{a}_{s'}^2 \to \Q \to 0.$$

Chacun de ces systèmes locaux est sous-jacent à une variation de structures de Hodge sur $A_{s'}$.
En effet, sur $\Q[\pi_1(\overline{A}_{s'},e(s'))] / \mathfrak{a}_{s'}^2 $ les filtrations proviennent
de la th\'eorie des int\'egrales it\'er\'ees de Chen (voir la partie \ref{deflogabsolu})
et l'augmentation $\Q[\pi_1(\overline{A}_{s'},e(s'))] / \mathfrak{a}_{s'}^2 \to \Q$
est sous-jacente \`a un morphisme de variations de structures de Hodge de but
la variation triviale $\Q(0)$.
Les filtrations sur $\mathfrak{a}_{s'} / \mathfrak{a}_{s'}^2$ sont celles induites par
celles de $\Q[\pi_1(\overline{A}_{s'},e(s'))] / \mathfrak{a}_{s'}^2$. En fait,
$\mathfrak{a}_{s'} / \mathfrak{a}_{s'}^2$ est la variation constante sur $A_{s'}$
associ\'ee \`a $H_1(\overline{A_{s'}},\Q)$. \\

Ainsi, sur
chacune des fibres de $\overline{\pi}$, on dispose de filtrations pour les trois systèmes locaux.
Il existe trois variations de $\Q$-structure de Hodge admissibles
sur $A$ dont les systèmes locaux sous-jacents et les filtrations fibre à fibre coïncident avec les données
précédentes. \\

\begin{itemize}
\item[a)] On note (abusivement) $\mathfrak{a}_s / \mathfrak{a}_s^2$ le système local sur $\overline{B}$
          associé  à $\mathfrak{a}_s / \mathfrak{a}_s^2$
          muni de l'action de $\pi_1(\overline{S},s)$ par multiplication.
          L'isomorphisme canonique $\mathfrak{a}_s / \mathfrak{a}_s^2 \to H_1(\overline{A}_s,\Q)$
          fournit un isomorphisme de sytèmes locaux sur $\overline{S}$ entre
          $\mathfrak{a}_s / \mathfrak{a}_s^2$ et
          $\overline{(R^1 \pi_* \Q)^{\vee}}$. On rappelle que $\hh$ d\'esigne la variation de structures de Hodge
	  pures de poids $-1$ $(R^1 \pi_* \Q)^{\vee}$. On \'equipe le
	  $\Q$-syst\`eme local $\mathfrak{a}_s / \mathfrak{a}_s^2$ sur $\overline{A}$ de la structure
	  de $\Q$-variations de $\pi^* \hh$ ($\pi_1(\overline{A}_s,a)$
	  agit trivialement sur $\mathfrak{a}_s / \mathfrak{a}_s^2$).
\item[b)] Pour $\Q[\pi_1(\overline{A}_s,a)]/ \mathfrak{a}_s^2$, on est dans la situation
          d'une variation sur un espace de chemins (cf. [HZ] et la construction de $\Log$ \cite[I-Thm 3.3]{w}).
          Les filtrations d\'efinies pr\'ec\'edemment fibre \`a fibre d\'efinissent donc  une $\Q$-variation admissible sur $A$
          not\'ee 
          $\Log_{A,s}^{(1)}$ ou  simplement  $\Log_{s}$ lorsque le sch\'ema ab\'elien est implicite,
          et on munit
	  le $\Q$-syst\`eme local $\Q[\pi_1(\overline{A}_s,a)]/ \mathfrak{a}_s^2$ sur $\overline{A}$
          de cette structure.
\item[c)] Pour $\Q$, on choisit $\Q(0)$. \\
\end{itemize}

Les morphismes figurant dans la suite exacte $(\Sigma)$ respectent les filtrations, d'où une suite exacte
dans $VSHM(A)$:
$$(\Sigma') \quad \quad \quad  0 \to  \pi^* \hh  \to \Log_{s}^{(1)}
 \to \Q(0)\to 0.$$
On remarque que, puisque $\text{Hom}_{VSHM(A)}( \Q(0), \pi^* \hh) = 0$, en raison des
poids, le terme m\'edian d'une suite exacte courte repr\'esentant un \'el\'ement de
$\mbox{Ext}^1_{VSHM(A)} (\Q(0) , \pi^* \hh )$ dans la description des Ext-groupes de Yoneda est bien d\'efini
\`a isomorphisme unique pr\`es. On notera ainsi \'egalement $\Log^1_s$ la classe dans
 $\mbox{Ext}^1_{VSHM(A)} (\Q(0) , \pi^* \hh )$ de la suite exacte $(\Sigma')$. \\

On cherche maintenant \`a caract\'eriser $\Log^{(1)}_s$ dans le groupe d'extensions
$\mbox{Ext}^1_{VSHM(A)} (\Q(0) , \pi^* \hh )$. La suite spectrale de Leray de la composition
$\text{RHom}_{MHM(S)}( \Q(0), \cdot ) \circ \pi_*$
appliquée à $\pi^* \hh$ donne la suite exacte courte scindée:
 $$ \xymatrix{ 0 \ar[r] &
           \text{Ext}^1_{MHM(S)}( \Q(0), \hh )
                     \ar[r]_-{\pi^*} &
  H^1 \text{RHom}_{MHM(A)}(  \Q(0) ,  \pi^* \hh )
  \ar[d] \ar@/_/[l]_-{e^*} & & \\
  & & Hom_{MHM(S)} (  \Q(0) , H^1 \pi_*  \pi^*  \hh    )
      \ar[r] & 0.
}$$

En effet pour des raisons de poids, $\text{Ext}^2_{MHM(S)}( \Q(0), \hh ) = 0$.
D'autre part $ H^1 \pi_*  \pi^*  \hh  =  \hh \otimes \hh^{\vee}$ (formule de projection)
et donc, par dualit\'e,
$$ Hom_{MHM(S)} (  \Q(0) , H^1 \pi_*  \pi^*  \hh    ) =
   Hom_{MHM(S)} ( \hh , \hh ) = \text{End}_{VSHM(S)}(\hh). $$
De plus, le foncteur exact et pleinement fid\`ele canonique $\iota_S \colon VSHM(S) \to MHM(S)$
induit un isomorphisme entre $\text{Ext}^1_{VSHM(S)}( \Q(0), \hh )$ et
$\text{Ext}^1_{MHM(S)}( \Q(0), \hh )$.
C'est une cons\'equence de la remarque suivant le Th\'eor\`eme 3.27 de \cite{s}.
De fa\c{c}on analogue, le foncteur $\iota_A \colon VSHM(A) \to MHM(A)$ induit
un isomorphisme entre $\text{Ext}^1_{VSHM(A)}( \Q(0), \pi^*\hh )$ et
$\text{Ext}^1_{MHM(A)}( \Q(0), \pi^* \hh ) .$
La suite exacte pr\'ec\'edente se r\'e\'ecrit donc comme suit:

$$ \xymatrix{ 0 \ar[r] &
           \text{Ext}^1_{VSHM(S)}( \Q(0), \hh )
                     \ar[r]_-{\pi^*} &
 \text{Ext}^1_{VSHM(A)}(  \Q(0) ,  \pi^* \hh )
  \ar[d]^{\partial} \ar@/_/[l]_-{e^*} & & \\
  & & \text{End}_{VSHM(S)}(\hh)
      \ar[r] & 0.
}$$

\begin{prop} $-$ \label{caraclog1}
L'extension $\Log^{(1)}_s$ v\'erifie $e^*\Log^{(1)}_s = 0$ et $\partial \Log^{(1)}_s = \text{Id}_{\hh}$.
\end{prop}

\begin{dem} \\

\begin{itemize}
\item[a)]
On consid\`ere la suite exacte $(\Sigma)$ comme suite exacte de $\pi(\overline{S},s)$-modules.
On remarque que le morphisme $1^{(1)}_s$ est $\pi(\overline{S},s)$-\'equivariant. Ainsi, $1^{(1)}_s$
fournit un scindage de la suite exacte $ \overline{e}^* (\Sigma)$ au niveau des $\Q$-syst\`emes locaux.
Pour d\'emontrer que $e^*\Log^{(1)}_s = 0$ est nul, il suffit donc de voir que $1^{(1)}_s$ respecte les filtrations
ce qui peut se v\'erifier fibre \`a fibre. Si $s' \in S(\C)$, le choix d'un chemin de $s'$ \`a $s$ fournit
une identification de $\Log^{(1)}_s$ et $\Log^{(1)}_{s'}$. Aussi suffit-il de montrer que $1^{(1)}_s$ 
respecte
les filtrations de la fibre en $s$. En utilisant la fonctorialit\'e des constructions pr\'ec\'edentes,
cette assertion est cons\'equence du fait que le morphisme
$$ 1 \colon \Q(0) \to \Q[\pi(\overline{A_s},e(s))]^{\; \widehat{\;}}$$
est un morphisme dans la cat\'egorie $SHM$ (cf. partie \ref{deflogabsolu}).\\

\item[b)] Pour d\'emontrer que $\partial \Log^{(1)}_s = \text{Id}_{\hh}$, il suffit de prouver cette
          identit\'e fibre \`a fibre. Comme en a), on r\'eduit ainsi l'\'etude au cas o\`u $A$ est une
	  vari\'et\'e ab\'elienne.\\
	  Dans ce cas $\Log^{(1)}$ correspond \`a la suite exacte
	  $$ 0 \to \mathfrak{a} / \mathfrak{a}^2 \to \Q[\pi_1(\overline{A},e)]  /  \mathfrak{a}^2 \to \Q(0) \to 0,$$
          o\`u $\mathfrak{a}$ est le noyau de l'augmentation
	  $\varepsilon \colon \Q[\pi_1(\overline{A},e)]  \to \Q$.
          Le morphisme $\partial\Log^{(1)}$
	  se d\'eduit par dualit\'e d'un morphisme de bord $\delta$
	  apparaissant dans la suite exacte longue de cohomologie associ\'ee au triangle
	  distingu\'e
          $\pi_* \pi^*  \mathfrak{a} / \mathfrak{a}^2   \to \pi_*  \Q[\pi_1(\overline{A},e)]  /  \mathfrak{a}^2
	   \to \pi_*\Q(0) \to
           \pi_* \pi^* \mathfrak{a} / \mathfrak{a}^2 [1] \quad :$
          $$\xymatrix{0 \to H^0 \pi_* \pi^* \mathfrak{a} / \mathfrak{a}^2  \to
          H^0  \pi_*   \Q[\pi_1(\overline{A},e)]  /  \mathfrak{a}^2   \to
             H^0  \pi_* \pi^* \Q(0)   \ar[r]^-{\delta} &
    H^1  \pi_* \pi^*  \mathfrak{a} / \mathfrak{a}^2 \to \dots }\quad .$$
    On souhaite donc d\'emontrer que $For(\delta)$ co\"incide avec $\text{Id}_{\mathfrak{a} / \mathfrak{a}^2 }$
    via l'identification

   $$ \text{Hom}_{\Q\text{-vect}} ( \Q , H^1 \overline{\pi}_* \overline{\pi}^* \mathfrak{a} / \mathfrak{a}^2  )
      =   \text{Hom}_{\Q\text{-vect}} ( \mathfrak{a} / \mathfrak{a}^2 , \mathfrak{a} / \mathfrak{a}^2).$$

  La compatibilité, via le foncteur $For$, entre les
  formalismes des six foncteurs au niveau des modules de Hodges d'une part, et au niveau topologique d'autre part,
  implique
  que $For(\delta)$ apparaît dans la suite exacte longue de cohomologie
  associée au triangle distingué
  $R \overline{\pi}_* \mathfrak{a} / \mathfrak{a}^2  \to R \overline{\pi}_*
    \Q[\pi_1(\overline{A},e)]  /  \mathfrak{a}^2  \to
      R \overline{\pi}_* (\Q) \to R \overline{\pi}_*
     (\mathfrak{a} / \mathfrak{a}^2) [1] \quad :$
    $$\xymatrix{0 \to H^0 R\overline{\pi}_* \mathfrak{a} / \mathfrak{a}^2 \to
                   H^0 R\overline{\pi}_*  \Q[\pi_1(\overline{A},e)]  /  \mathfrak{a}^2     \to
          H^0 R\overline{\pi}_*  \Q \ar[rr]^-{For(\delta)}  & &
       H^1  R\overline{\pi}_*  \mathfrak{a} / \mathfrak{a}^2 \to \dots }\quad .$$

Comme $\overline{A}$ est un tore, c'est un $K(\Gamma,1)$. On peut donc utiliser la cohomologie
  du groupe $\pi_1(\overline{A},e)$  pour calculer $For(\delta)$ .
  Ce dernier est présent dans la suite exacte longue de cohomologie associée
  à la suite exacte courte de représentations de $\pi_1(\overline{A},e) \quad $
  $ 0 \to \mathfrak{a} / \mathfrak{a}^2 \to  \Q[\pi_1(\overline{A},e)]  /
\mathfrak{a}^2  \to \Q \to 0 \quad :$
  $$ \xymatrix{0 
\to \mathfrak{a} / \mathfrak{a}^2 \to 
(\Q[\pi_1(\overline{A},e)]  /  \mathfrak{a}^2)^{\pi_1(\overline{A},e)}
  \to \Q
  \ar[rr]^-{For(\delta)}  &    & 
H^1(\pi_1(\overline{A},e), \mathfrak{a}/\mathfrak{a}^2)  \to \dots} \quad.$$

Pour calculer $For(\delta)$, on introduit le diagramme suivant:
$$ \xymatrix{
        0 \ar[r]      &  \mathfrak{a} /  \mathfrak{a}^2 \ar[r]^-{\iota}  \ar[d]^-{d'^{0}} &
             \Q[\pi_1(\overline{A},e)] / \mathfrak{a}^2               
\ar[r]^-{\overline{\epsilon^{(1)}}}   \ar[d]^-{d^{0}}
 & \Q \ar[d]^{d^{\prime\prime0}} \ar[r]
& 0 \\
           & L_1(\mathfrak{a} /  \mathfrak{a}^2 ) \ar[r]^-{\iota_*}  \ar[d]^-{d'^{1}} 
&             L_1(\Q[\pi_1(\overline{A},e)] / \mathfrak{a}^2 ) 
\ar[r]^-{(\overline{\epsilon^{(1)}})_*}  \ar[d]^-{d^1} &
             L_1(\Q)    \ar[d]^{d^{\prime\prime1}} &  \\
              & L_2(\mathfrak{a} /  \mathfrak{a}^2 ) \ar[r]^-{\iota_*} 
&
             L_2(\Q[\pi_1(\overline{A},e)] / \mathfrak{a}^2) 
\ar[r]^-{(\overline{\epsilon^{(1)}})_*}
       &
             L_2(\Q)  & \\
             \\
         }    \\$$

\noindent dans lequel $L_i(V)$ ($i \in \N$, $V$ $\Q$-vectoriel) 
désigne l'ensemble des applications du produit $\pi_1(\overline{A},e)^i$
dans $V$ (muni de la structure de $\Q$-vectoriel \'evidente) et o\`u les morphismes
verticaux sont les diff\'erentielles usuelles.  On se donne de plus un isomorphisme
de groupes ab\'eliens
$\pi_1(\overline{A},e) \simeq \Z^{2d}$ et on identifie $\Q[\pi_1(\overline{A},e)] $ 
\`a  $\Q[X_1,..,X_{2d},X_1^{-1},..,X^{-1}_{2d}]$. On v\'erifie alors que le morphisme suivant
 est un isomorphisme de groupes ab\'eliens.
$$ \gamma \colon \pi_1(\overline{A},e) \to \mathfrak{a}/\mathfrak{a}^2, \quad
     (n_1,\dots,n_{2d}) \mapsto n_1 \overline{ (X_1 - 1 )}  + \dots + n_{2d} \overline{(X_{2d} - 1 )} $$

On a:
$$ \begin{array}{llll}
   H^1(\pi_1(\overline{A},e), \mathfrak{a}/\mathfrak{a}^2) & = &  
   \mbox{Ker}(d'^{1}) / \mbox{Im}(d'^0) & \\ 
   & = & \text{Hom}_{\Z\text{-Mod}}( \pi_1(\overline{A},e) , 
\mathfrak{a}/\mathfrak{a}^2 ) & (\pi_1(\overline{A},e) \text{ agit tivialement sur } 
 \mathfrak{a}/\mathfrak{a}^2) \\    
  & \underset{(*)}{\simeq} & \text{Hom}_{\Z\text{-Mod}}( \mathfrak{a}/\mathfrak{a}^2 , 
\mathfrak{a}/\mathfrak{a}^2 ) & (\text{via l'isomorphisme } \gamma )
\end{array} $$
et pour tout $x \in \Q$, 
$For(\delta)(x) \in H^1(\pi_1(\overline{A},e), \mathfrak{a}/\mathfrak{a}^2)$
est donn\'e par la classe dans 
$H^1(\pi_1(\overline{A},e), \mathfrak{a}/\mathfrak{a}^2)$ de l'\'el\'ement 
(bien d\'efini) 
$\iota_*^{-1} \; d^0 \; (\overline{\epsilon^{(1)}})^{-1} \; (x) \in \Ker(d'^1).$\\
En raison de la $\Q$-lin\'earit\'e, il suffit de montrer que 
$For(\delta)(1)$ co\"incide avec 
$Id_{\mathfrak{a}/\mathfrak{a}^2}$ via l'identification $(*)$. Or un calcul \'el\'ementaire montre que
$ \iota_*^{-1} \; d^0 \; (\overline{\epsilon^{(1)}})^{-1} \; (1) $ est l'isomorphisme $\gamma$.
\end{itemize}

\end{dem}

On note que les deux propri\'et\'es pr\'ec\'edentes caract\'erisent l'extension/la $\Q$-variation
$\Log^{(1)}_s$. D'autre part, $\Log^{(1)}_s$ est \'equip\'e d'un morphisme canonique  $\varepsilon^{(1)}_s \colon \Log^{(1)}_s \to \Q(0)$ dans la cat\'egorie
$VSHM(A)$ et d'un morphisme 
$1^{(1)}_s \colon \Q(0) \to e^* \Log^{(1)}_s$ dans $VSHM(S)$ .
Comme on l'a remarqu\'e pr\'ec\'edemment, le couple
$(\Log^{(1)}_s , \varepsilon^{(1)}_s )$ est rigide et par suite ne d\'epend pas du
choix de $s$. On s'autorisera donc \`a noter simplement $(\Log^{(1)}, \varepsilon^{(1)})$ le couple $(\Log^{(1)}_s , \varepsilon^{(1)}_s )$.

\begin{notas} $-$ \\

\begin{tabular}{cp{14cm}}
$\Log^{(n)}$ & $:= \text{Sym}^n \Log^{(1)}$, pour $n \in \N$, \\
$c_n(\chi)$  & l'application de $\Sym^n \; V$ vers $\Sym^{n-1} \; V$, qui associe \`a
               $[v_1 \otimes \dots \otimes v_n]$ l'\'el\'ement
              $$\frac{1}{n!} \; \sum_{\sigma \in \mathfrak{S}_n}
                \; \chi(v_{\sigma(1)}) \;  [v_{\sigma(2)} \otimes .. \otimes v_\sigma(n)],$$
              pour $n \in \N^*$, o\`u  $V$ est une $\Q$-repr\'esentation de dimension finie de
	      $\pi_1(A,a)$ et $\chi: V \to \Q$ une forme linéaire $\pi_1(A,a)$-invariante.

\end{tabular}
\end{notas}

On considère le pro-objet de $VSHM(A)$
$ \displaystyle \lim_{\overset{\longleftarrow}{n \geq 0}} \Log^{(n)}$
dont les morphismes de transition sont donnés au niveau des $\Q$-syst\`emes locaux
par les $c_n(\overline{\varepsilon^{(1)}})$,  $n \geq 0$.

\begin{rema} $-$ Dans \cite{k}, Kings d\'efinit le logarithme du sch\'ema ab\'elien $A/S$ comme
\'etant $ \displaystyle \lim_{\overset{\longleftarrow}{n \geq 0}} \Log^{(n)}$.
\end{rema}

On d\'emontre maintenant
que $\displaystyle \lim_{\overset{\longleftarrow}{n \geq 0}} \Log^{(n)}$ est isomorphe \`a $\Log$ d\'efini
dans la partie \ref{deflogrel}.
Pour tout $n \geq 1$, soit $1_s^{(n)}: \Q(0) \to e^* \Log_s^{(n)}$ le morphisme induit par:

$$ \begin{array}{ccccl}
      \overset{n}{\otimes}  \; 1^{(1)}_s: &
       \overset{n}{\otimes} \; \Q(0) & \to & \overset{n}{\otimes} \; e^* \Log_s^{(1)}
     \end{array}.
          $$

On applique alors le Théorème \ref{thmutile} pour associer à $1_s^{(n)}$, pour $n \in \N$, le morphisme
$\varphi_s^{(n)}$
$$ \varphi_s^{(n)}  : \Log \to \Log_s^{(n)}.$$
Puisque pour tout $n\geq1$,
$c_n(\overline{\varepsilon^{(1)}}) \circ \overline{1_s^{(n)}}  = \overline{1_s^{(n-1)}} $,
$\left(\varphi_s^{(n)}\right)_{n \in \N}$ définit un morphisme
dans pro-$VSHM(A)$
$$\displaystyle \varphi_s : \Log\to \lim_{\overset{\longleftarrow}{n \geq 0}}  \Log^{(n)}.$$

\begin{prop} $-$ Le morphisme $\varphi_s$ est un isomorphisme.
\end{prop}

\begin{dem}
Il suffit de prouver que $\varphi_s$ induit un morphisme sur chacune des fibres. Par fonctorialit\'e de
la construction de $\varphi_s$, on r\'eduit l'assertion au cas o\`u $A$ est une vari\'et\'e ab\'elienne.
Dans ce cas, on supprime l'indice $s$ dans les notations.
Il est suffisant de prouver que
$\overline{\varphi} : \Q[\pi_1(\overline{A},e)]^{\; \widehat{\;}}
\to \displaystyle  \lim_{\overset{\longleftarrow}{n \geq 0}}
\Sym^n (\Q[\pi_1(\overline{A},e)] / \mathfrak{a}^2)$
est un isomorphisme de $\Q$-vectoriels, o\`u $\mathfrak{a}$ d\'esigne  le noyau de l'augmentation
de $\Q[\pi_1(\overline{A},e)]$.
On fixe un isomorphisme $\pi_1(\overline{A},e) \simeq \Z^{2d}$. Celui-ci détermine un isomorphisme
 $\Q[\pi_1(\overline{A},e)] \simeq \Q[X_1,X_1^{-1}, ..,X_{2d}, X_{2d}^{-1}].$
 Soit $n \in \N$.
Puisque $\Log_A^{(n)}$ est $(n+1)$-unipotente (sa filtration par le poids a $(n+1)$-gradu\'es non triviaux qui sont
des variations constantes),
le morphisme $\varphi^{(n)}$ se factorise donc à travers la projection
$\Q[\pi_1(\overline{A},e)]^{\; \widehat{\;}} \to \Q[\pi_1(\overline{A},e)] / \mathfrak{a}^{n+1}$.
L'action de $\pi_1(\overline{A},e)$ sur $ \Q[\pi_1(\overline{A},e)] / \mathfrak{a}^2$
étant donnée par la multiplication, on en déduit que $\overline{\varphi^{(n)}}$ est donn\'e par la composition:

   $$ \begin{array}{ccccc}
       \Q[\pi_1(\overline{A},e)]^{\; \widehat{\;}} &  \to & \Q[\pi_1(\overline{A},e)]  / \mathfrak{a}^{n+1}  &
       \overset{\psi^{(n)}}{\to} &  \Sym^n (\Q[\pi_1(\overline{A},e)] / \mathfrak{a}^2) \\
          & & X_1^{i_1}..X_{2d}^{i_{2d}} & \mapsto &
        [ [X_1^{i_1}..X_{2d}^{i_{2d}}] \otimes .. \otimes [ [X_1^{i_1}..X_{2d}^{i_{2d}}] ] \end{array}.$$
On remarque que $\psi^{(0)}= \text{Id}_{\Q}$. On considère le diagramme commutatif suivant:
$$\xymatrix{
0 \ar[r] & \mathfrak{a}^{n+1} / \mathfrak{a}^{n+2} \ar[r] \ar[d]^-{\psi^{(n+1)}_{|}}   & 
           \Q[\pi_1(\overline{A},e)]  / \mathfrak{a}^{n+2}   \ar[r] \ar[d]^-{\psi^{(n+1)}}
& \Q[\pi_1(\overline{A},e)]  / \mathfrak{a}^{n+1}    \ar[r] \ar[d]^-{\psi^{(n)}} & 0 \\
0 \ar[r] & \Sym^{n+1} ( \mathfrak{a} / \mathfrak{a}^2) \ar[r]_-{i_{n+1}} &
             \Sym^{n+1} (\Q[\pi_1(\overline{A},e)] / \mathfrak{a}^2) \ar[r]_-{
            \overset{\quad}{c_{n+1}(\overline{\varepsilon^{(1)}})}} &
\Sym^n (\Q[\pi_1] / \mathfrak{a}^2) \ar[r] & 0}$$
où le morphisme $i_{n+1}$ est induit par l'inclusion
$\mathfrak{a} / \mathfrak{a}^2 \subset \Q[\pi_1(\overline{A},e)] / \mathfrak{a}^2$. On prouve maintenant que $\psi^{(n+1)}_{|}$
est un isomorphisme. La famille
$$ \left( [(X_1-1)^{i_1}..(X_{2d}-1)^{i_{2d}}]\right)_{\{ (i_1,..,i_{2d}) \in \N^{2d} \; / \;  i_1+..+i_{2d}=n+1 \}}$$
est une base de $\mathfrak{a}^{n+1} / \mathfrak{a}^{n+2}$. Pour tout
 $(i_1,..,i_{2d}) \in \N^{2d} \; / \;  i_1+..+i_{2d}=n+1$, on a
$$ \psi^{(n+1)} ([(X_1-1)^{i_1}..(X_{2d}-1)^{i_{2d}}]) =
    [\overset{i_1}{\otimes} [X_1-1] \otimes .. \otimes  \overset{i_{2d}}{\otimes} [X_{2d}-1]].$$
Or $\{[X_1-1],..,[X_{2d}-1]\}$ est une famille libre de $\Q[\pi_1(\overline{A},e)] / \mathfrak{a}^2$. Donc $\psi^{(n+1)}_{|}$
est injective. On conclut à la bijectivité à l'aide des dimensions. À l'aide d'une récurrence, on déduit donc que les $\psi^{n}$
sont des isomorphismes.

\end{dem}

\subsection{Le pro-syst\`eme local sous-jacent au logarithme}

Dans cette partie on d\'ecrit le pro-syst\`eme local de $\R$-vectoriels
$\overline{\Log}_{\R}$ \`a l'aide d'un  pro-fibr\'e vectoriel \`a connexion int\'egrable
sur $A^{\infty}$ en utilisant la construction de Levin (cf. \cite[Part 2]{l}).

\subsubsection{Le fibr\'e tangent d'une famille de tores r\'eels}

Soit $B$ une vari\'et\'e diff\'erentielle.

\begin{defi} $-$ Une famille de groupes de Lie r\'eels  au-dessus de $B$ est la donn\'ee
d'une fibration en tores r\'eels $p : X \to B$ et de trois morphismes de vari\'et\'es diff\'erentielles
$0: B \to X$  (unité) section de $p$,
$m: X \times_B X \to X$  (multiplication) compatible avec les projections sur $B$,
$i: G \to G$ (inverse) tel que $p \circ i = p$,  de sorte que le quadruplet
$(p,0,m,i)$ d\'efinit un objet en groupes dans la cat\'egorie des vari\'et\'es diff\'erentielles
au-dessus de $B$. On note que $p$ \'etant une fibration, le produit fibr\'e $X \times_B X$ 
dans la cat\'egorie des vari\'et\'es diff\'erentielles est bien d\'efini.
Si les fibres de $p$ sont des tores r\'eels, on dit que $(p,0,m,i)$ est une famille de tores r\'eels au-dessus de
$B$.
\end{defi}

On a une notion \'evidente de morphisme entre familles de groupes de Lie r\'eels au-dessus de $B$.
On note $Lie_{/B}$ la cat\'egorie des familles de groupes de Lie r\'eels au-dessus de $B$. \\

Soit $\Gamma$ un syst\`eme local de groupes ab\'eliens libres de rang fini au-dessus de $B$.
Il existe une construction classique qui permet d'associer au faisceau de $\OO_B$-modules
localement libres $\Gamma \otimes \OO_B$ un fibr\'e vectoriel au-dessus de $B$ que l'on note
$E(\Gamma)$ dont la fibre au-dessus de $b \in B$ est $(\Gamma_b)_\R$ et tel que
le faisceau des sections de  $E(\Gamma)$ est $\Gamma \otimes \OO_B$.

\begin{fait} $-$
On peut, de mani\`ere analogue,
construire \`a partir de $\Gamma$ une famille de tores r\'eels au-dessus de $B$,
not\'ee $p \colon E(\Gamma)/\Gamma \to B$, dont la fibre au-dessus de $b \in B$
 est $(\Gamma_b)_\R / \Gamma_b$.
Par construction $E(\Gamma)/\Gamma$ se trouve être muni d'un morphisme canonique
$$q \colon E(\Gamma) \to E(\Gamma)/\Gamma $$ dans $Lie_{/B}$ qui est universel, i.e.
pour tout morphisme $r \colon E(\Gamma) \to X$ dans $Lie_{/B}$ tel que $\Gamma$ est un sous-faisceau du
faisceau des sections du noyau de $r$, il existe un unique morphisme $\overline{r} \colon E(\Gamma) / \Gamma \to X$
dans $Lie_{/B}$ tel que $\overline{r} \circ p = r$. De plus, on a une d\'ecomposition canonique du
fibr\'e tangent de $E(\Gamma) / \Gamma$:
$$ TE(\Gamma) / \Gamma = p^* E(\Gamma) \oplus p^* TB.$$
\end{fait}

En fait, toute famille de tores r\'eels au-dessus de $B$ est isomorphe \`a une famille
de tores r\'eels au-dessus de $B$ ainsi construite.
Soit $(p \colon X \to B,0,m,i)$ une famille de tores r\'eels au-dessus de $B$.
L'exponentielle fibre \`a fibre d\'efinit une application diff\'erentielle
$\exp \colon 0^* TX_{/B} \to X$ o\`u $TX_{/B}$ est le noyau de $Tp \colon TX \to TB$.
Le faisceau des sections du noyau de $\exp$ s'identifie \`a $\Gamma:=(R^1 p_* \Z)^{\vee}$. 
Le faisceau des sections de $0^* TX_{/B}$ est donc canoniquement isomorphe \`a  $\Gamma \otimes_{\Z} \OO_{B}$.
Ainsi, en factorisant par $q \colon E(\Gamma) \to E(\Gamma)/\Gamma$, on en d\'eduit
un morphisme $\overline{\exp} \colon TE(\Gamma) / \Gamma  \to X$ dans $Lie_{/B}$ qui
est un isomorphisme. En effet, c'est un isomorphisme sur les fibres au-dessus de $B$.
On obtient donc une d\'ecomposition canonique du fibr\'e tangent de $X$
$$ TX = p^* E(\Gamma) \oplus p^* TB.$$

\subsubsection{Description de $\overline{\Log}_{\R}$}

La vari\'et\'e diff\'erentielle $\pi^{\infty} \colon A^{\infty}$ au-dessus de $S^{\infty}$
est munie d'une structure de famille de tores r\'eels au-dessus de $S^{\infty}$ h\'erit\'ee
des lois de structure du sch\'ema ab\'elien $A/S$.
On note $\Gamma$ le syst\`eme local $(R^1 \overline{\pi} \Z)^{\vee}$ sur $S^{\infty}$
et simplement $E$ le fibr\'e vectoriel r\'eel $E(\Gamma)$ sur $S^{\infty}$.
On remarque qu'avec les notations introduites, on a  $\overline{\hh}= \Gamma_{\Q}$.\\

Le morphisme $\partial \colon \text{Ext}^1_{VSHM(A)}(  \Q(0) ,  \pi^* \hh ) \to
 \text{End}_{VSHM(S)}(\hh)$ de la Proposition \ref{caraclog1}
 a \'et\'e construit en consid\'erant la th\'eorie
des modules de Hodge mixtes. On peut, de fa\c{c}on analogue, construire un morphisme
$$ For(\partial)_{\R} \colon \text{Ext}^1_{\F_{\R}(A)} ( \R , \overline{\pi}^* \Gamma_{\R}) \to
   \text{End}_{\F_{\R}(S) }(\Gamma_{\R} )$$
en se pla\c{c}ant, cette fois, au niveau topologique. On d\'eduit de la Proposition
\ref{caraclog1} que le syst\`eme local de $\R$-vectoriels $(\overline{Log^{(1)}})_{\R}$
est caract\'eris\'e par
$$ \overline{e}^*(\overline{Log^{(1)}})_{\R} = 0 \text{ et }
   For(\partial)_{\R} ( (\overline{Log^{(1)}})_{\R} ) = \text{Id}_{\Gamma_{\R}}.$$

On a vu dans la partie pr\'ec\'edente que l'exponentielle fibre \`a fibre induisait une  d\'ecomposition
du fibr\'e tangent de $A^{\infty}$: $TA^{\infty} = (\pi^{\infty})^* E \oplus (\pi^{\infty})^* TB$. On note
$\nu$ la 1-forme diff\'erentielle sur $A^{\infty}$ \`a valeurs dans $(\pi^{\infty})^* E$ correspondant \`a la
projection canonique de $TA^{\infty}$ sur $(\pi^{\infty})^*E$.

\begin{lemme} \label{lemmenu} $-$ Soit $\nabla_{GM}$ la connexion de Gauss-Manin sur $E$.
La forme $\nu$ est ferm\'ee, i.e. $\nabla_{GM}(\nu)=0$,
et sa classe $[\nu]$ dans $H^1(\overline{A}, \overline{\pi}^* \Gamma_{\R})$ v\'erifie
 $For(\partial)_{\R} ([\nu]) = \text{Id}_{\Gamma_{\R}}$.
\end{lemme}

\begin{dem} \\

\begin{itemize}
\item[a)] Pour la preuve de la premi\`ere assertion, on renvoie le lecteur \`a \cite[p. 216]{l}. \\

\item[b)] Pour la seconde, il suffit de v\'erifier l'identit\'e sur les fibres. Les  constructions
\'etant fonctorielles, on s'est ainsi ramen\'e \`a prouver la relation dans le cas o\`u
$A$ est une vari\'et\'e ab\'elienne.
Dans ce cas, le morphisme $For(\partial)_{\R}$ est donn\'e par la composition:
$$ \begin{array}{rcl}  H^1( \overline{A}, H_1(\overline{A},\R)) =   H^1( \overline{A} , \R) \otimes H_1(\overline{A} , \R) & \to & \text{End}( H_1(\overline{A} , \R)). \\
                                                                        \omega \otimes c & \mapsto &  ( c' \mapsto <\omega,c'> c) \end{array}$$
On fixe un isomorphisme $A^{\infty} \simeq \R^{2d} / \Z^{2d}$. On obtient alors des coordonn\'ees et on exprime $\nu$ relativement \`a celles-ci. On calcule
$For(\partial)_{\R}([\nu])$ \`a l'aide de la composition donn\'ee ci-dessus pour \'etablir 
$For(\partial)_{\R}([\nu] = \text{Id}_{ H_1(\overline{A} , \R)}$.
\end{itemize}
\end{dem}

On peut maintenant expliciter un fibré à connexion candidat pour représenter le système local
$(\overline{\Log^{(1)}})_{\R}$.
On considère le fibré vectoriel $E' := \OO_{A^{\infty}} \oplus (\pi^{\infty})^* E $
muni de la connexion $\nabla^1$:
$$ \begin{array}{rccc}
    \nabla^1 \; : \; &  \OO_{A^{\infty}} \oplus  (\pi^{\infty})^* E & 
                 \longrightarrow  &
                     \Omega^1_{A^{\infty}} \oplus  \Omega^1_{A^{\infty}} \otimes  (\pi^{\infty})^* E .\\
    & (f, g \otimes h) & \mapsto & (df , dg \otimes h + f \nu) \end{array}$$
La connexion $\nabla^1$ est plate ($\nu$ est fermée). Le faisceau $\mathbb{E} := \mbox{Ker}(\nabla^1)$
est donc un système local.
On a une suite exacte de fibrés vectoriels munis de connexions:
$$\begin{array}{ccccccccc}
0 & \to & ((\pi^{\infty})^* E , \nabla_{GM}) & \to &
\left( \OO_{A^{\infty}} \oplus (\pi^{\infty})^*E , \nabla^1 \right) & \to &
(\OO_{A^{\infty}} , d) & \to & 0 ,\\
 & & g \otimes h & \mapsto & (0,g \otimes h) & & & &  \\
             & &  & &            ( f , g \otimes h) & \mapsto & f & & \\ \end{array} $$
où $\nabla_{GM}$ désigne la connexion de Gauss-Manin.
 Celle-ci correspond à une suite exacte de systèmes locaux
$ 0 \to \overline{\pi}^* \overline{\hh} \to \mathbb{E} \to \R \to 0$
dont la classe dans
$\mbox{Ext}^1_{\F_{\R}(A)}(\R,
 \overline{\pi}^* \Gamma_{\R})$ est notée $[\mathbb{E}]$.

\begin{prop}{}  \label{desclog1top} On a les identités suivantes:
$$  \overline{e}^*([\mathbb{E}])=0 \quad \mbox{ et } \quad
    For(\partial)_{\R}([\mathbb{E}])= \text{Id}_{\Gamma_{\R}}, $$
i.e. $\mathbb{E} = (\overline{\Log^{(1)}})_{\R}$.
\end{prop}

\begin{dem} La premi\`ere identit\'e est \'evidente. Pour d\'emontrer la deuxi\`eme,
on utilise la résolution de $\mathbb{E}$ construite à partir
de $(E',\nabla^1)$ pour expliciter $[\mathbb{E}] \in H^1(\overline{A} , \overline{\pi}^* \Gamma_{\R})$:
$$\xymatrix{
& \overline{\pi}^* \overline{\hh} \ar[r] \ar[d] &  \mathbb{E}  \ar[r]  \ar[d] &
 \R \ar[d]^-i \\
& (\pi^{\infty})^* E  \ar[r] \ar[d] &
\OO_{A^{\infty}} \oplus  (\pi^{\infty})^* E \ar[r]^-p  \ar[d]^-{\nabla^1} &
 \OO_{A^{\infty}} \ar[r] \ar[d] & 0 .\\
0 \ar[r] &  \Omega^1_{A^{\infty}} \otimes (\pi^{\infty})^* E \ar[r]^-j \ar[d]^-{d^1} &
\Omega^1_{A^{\infty}} \oplus \left( \Omega^1_{A^{\infty}} \otimes  (\pi^{\infty})^* E \right) \ar[r] \ar[d] &
 \Omega^1_{A^{\infty}} \ar[d] & \\
&\Omega^2_{A^{\infty}} \otimes (\pi^{\infty})^*  E   \ar[r] &
\Omega^2_{A^{\infty}} \oplus  \left( \Omega^2_{A^{\infty}} \otimes  (\pi^{\infty})^* E  \right) \ar[r]  &
\Omega^2_{A^{\infty}}  & \\ }$$
Alors $j^{-1} \;  \nabla^1 \; p^{-1} \; i \; (1)$ est dans $\mbox{Ker}(d^1)$ et sa classe dans
  $H^1(\overline{A} , \overline{\pi}^* \Gamma_{\R})$ coïncide
avec $[\mathbb{E}] \in H^1(\overline{A} , \overline{\pi}^* \Gamma_{\R})$.
Or $j^{-1} \;  \nabla^1 \; p^{-1} \; i \; (1)=\nu$. On conclut \`a l'aide
 du Lemme
\ref{lemmenu}.\\

\end{dem}

Apr\`es avoir obtenu cette description de $(\overline{\Log^{(1)}})_{\R}$, on \'etudie
$\overline{\Log}_{\R}$. Tout d'abord, $\nabla^1$ sur $E'$ induit une connexion
$\nabla^n$ sur $\Sym^n E'$, pour $n \in \N^*$.\\

Soit $\nu_n \colon \Sym^n (\pi^{\infty})^* E \to
                (\Sym^{n+1} (\pi^{\infty})^* E ) \otimes \Omega^1_{A^{\infty}}$, $n \geq 0$,
définie comme étant la composée:
$$\xymatrix{
     \Sym^n (\pi^{\infty})^* E \ar[r]^-{Id \otimes \nu} &
     ( \Sym^n (\pi^{\infty})^* E ) \otimes (\pi^{\infty})^* E \otimes \Omega^1_{A^{\infty}}
     \ar[r]^-{\underset{\quad}{mult \otimes Id}}
    & ( \Sym^{n+1} (\pi^{\infty})^*  E ) \otimes \Omega^1_{A^{\infty}}}.$$
On introduit alors le pro-fibré à connexion
$$ (\mathcal{G},\nabla) :=
    \left(  \prod\limits_{ n \geq 0 } \Sym^n (\pi^{\infty})^* E , \prod\limits_{ n \geq 0 }
            (\nabla_{GM}^n + \nu_n) \right),$$
o\`u  $\nabla_{GM}^n$ est la connexion de Gauss-Manin sur $\Sym^n (\pi^{\infty})^* E$.\

Soit $l \in \N^{\geq 2}$. On remarque que le sous-fibr\'e vectoriel
$ \displaystyle W_l := \prod\limits_{ k \geq l+1} \Sym^k (\pi^{\infty})^* E$ est stable par $\nabla$ et on d\'efinit:\\

\begin{tabular}{cl}
$(\G_l , \nabla_l)$ & le fibr\'e \`a connexion plate quotient $(\G, \nabla) / W_{l}$,\\
$p_l$ &  la projection canonique $(\G, \nabla) \to (\G_l , \nabla_l)$, \\
$p_{l+1,l}$ & la projection canonique  $(\G_{l+1}, \nabla_{l+1}) \to (\G_l , \nabla_l)$. \\ \\
\end{tabular}

On remarque que
les morphismes $p_l$ induisent
un morphisme de pro-fibr\'es vectoriels \`a connexions
plates
$$  p \colon (\G,\nabla) \to  \lim_{\overset{\longleftarrow}{l \geq 1}}  \;  ( \G_l , \nabla_l )$$
qui est un isomorphisme (les morphismes de transition de l'objet de droite sont les morphismes $p_{l+1,l}$).

\begin{prop} $-$ Il existe une famille d'isomorphismes de fibr\'es vectoriels \`a connexions plates
   $(\theta_n \colon (\G_n , \nabla_n) \to (\Sym^n E' , \nabla^n) )_{n \in \N}$ qui induit un isomorphisme
   de pro-fibr\'es vectoriels \`a connexions:
   $$ \theta \colon (\mathcal{G},\nabla) = \lim_{\overset{\longleftarrow}{n \geq 1}}  \;  ( \G_n , \nabla_n ) \to
                    \lim_{\overset{\longleftarrow}{n \geq 1}}  \;   
		   (\Sym^n  E'  , \nabla^n ). $$
   Et donc, le noyau de $\nabla$ s'identifie \`a $\overline{\Log}_{\R}$.
\end{prop}

\begin{dem}
On commence par remarquer que la derni\`ere assertion se d\'eduit de l'existence d'un tel isomorphisme
   $\theta$ et de
 la Proposition \ref{desclog1top}.\\
D'une part,
$\mathcal{G}_n = \underset{0 \leq k \leq n}{\bigoplus} \; \Sym^k (\pi^{\infty})^* E$
et d'autre part, on a un isomorphisme naturel de fibr\'es vectoriels:
$$ \begin{array}{lccc}
   \psi_n \; : & \underset{0 \leq k \leq n}{\bigoplus} \; Sym^k(\pi^{\infty})^* E &
               \to & Sym^k (\OO_{A^{\infty}} \oplus (\pi^{\infty})^* E ).\\
        & [h_1 \otimes .. \otimes h_k ] & \mapsto &
          [ 1 \otimes .. \otimes 1 \otimes h_1 \otimes .. \otimes h_k] \end{array}$$

Pour $n \geq 2$, $\psi_n$ n'est ni compatible avec les morphismes de transition, ni compatible avec les connexions.
 On corrige
ce défaut à l'aide d'un automorphisme $\alpha_n$ de
$ \underset{0 \leq k \leq n}{\bigoplus} \; \Sym^k (\overline{\pi^*\hh} \otimes \OO_{A^{\infty}})$ défini facteur
par facteur par une homothétie de rapport
   $$ \alpha_n^k :=  \frac{n!}{(n-k)!} , \quad n \in \N, \; 0 \leq k \leq n.$$ non nul.
Si on pose maintenant pour tout $n \geq 0$, $\theta_n:=\psi^n \circ \alpha_n$, on v\'erifie que
la famille $\left(  \theta_n \right)_{n \geq 0}$ est bien une famille d'isomorphismes de fibr\'es vectoriels
\`a connexions compatibles avec les morphismes de transitions. Ainsi, elle induit un isomorphisme
 $$ \theta \colon (\mathcal{G},\nabla) = \lim_{\overset{\longleftarrow}{n \geq 1}}  \;  ( \G_n , \nabla_n ) \to
                    \lim_{\overset{\longleftarrow}{n \geq 1}}  \;  ( \Sym^n E'  , \nabla^n ). $$

\end{dem}

\subsection{Propri\'et\'es du logarithme d'un sch\'ema ab\'elien} \label{proplog}

On rappelle que $\hh$ d\'esigne $(R^1 \pi_* \Q)^{\vee} \in Ob(VSHM(A))$.

\subsubsection{Gradu\'es par le poids}

Puisque l'on dispose d'une suite exacte canonique
$$ 0 \to \pi^* \hh \to \Log^{(1)} \overset{\varepsilon^{(1)}}{\to} \Q(0) \to 0,$$
on a une identification naturelle entre  le gradu\'e par le poids de
$\Log^{(1)}$ est $\Q(0) \oplus \pi^*\hh$. De cette propri\'et\'e et de l'isomorphisme
$\Log = \displaystyle \lim_{\overset{\longleftarrow}{n \geq 1}} \Sym^n \Log^{(1)}$,
o\`u les morphismes de transition dans le membre de droite sont induits par $\varepsilon^{(1)}$,
on d\'eduit que
$$ Gr^W \Log = \underset{n \geq 0}{\oplus} \; \Sym \; \pi^* \hh.$$

\subsubsection{Principe de scindage pour la section unit\'e} \label{secsc}

On a vu que $e^*\Log^{(1)} = \Q(0) \oplus \hh$ (cf. Proposition \ref{caraclog1}).
\`A nouveau en utilisant l'isomorphisme canonique
$\Log = \displaystyle \lim_{\overset{\longleftarrow}{n \geq 1}} \Sym^n \Log^{(1)}$,
on montre que
$$ e^* \Log = \prod\limits_{n \geq 0 } \; \Sym^n \hh.$$

\subsubsection{Principe de scindage pour une section de torsion} \label{principescindage}

Soit $x : S \to A$ une section de $N$-torsion.
Soit $[N] : A \to A$ l'isogénie donnée par la multiplication
par $N$. On applique $[N]^*$ à la suite exacte
$$ 0 \to \pi^* \hh \to \Log^{1} \to \Q(0) \to 0$$ pour obtenir
une suite exacte
$$ 0 \to \pi^* \hh =
     [N]^* \pi^* \hh \to [N]^* \Log^{1} \to \Q(0)  \to 0$$
dont on note  $[[N]^* \Log^{1}]$ la classe dans
$\text{Ext}_{VSHM(A)}^1( \Q(0) , \pi^* \hh)$.
De $e^* [\Log^{1}] = 0$ et $\partial[\Log^{1}] = \text{Id}_{\hh}$, on déduit
$e^* [[N]^* \Log^{1}] = 0$ et $\partial [[N]^* \Log^{1}] = \text{Id}_{\hh}$.
Par conséquent,
$[[N]^* \Log^{1}] = [\Log^{1}]$ (cf. Proposition \ref{caraclog1})
et donc $[N]^* \Log^{1} = \Log^{1}$. Ainsi
$$ [N]^* \Log =  \Log ( = \lim_{\overset{\longleftarrow}{n \geq 1}} \Sym^n \Log^{(1)}) . $$

\begin{prop}{\emph{[W, III-Prop 6.1]}} $-$ \label{scind}
$x^* \Log = \prod\limits_{k \geq 0} \Sym^k \hh $.
\end{prop}

\begin{dem} $ x^* \Log = x^* [N]^* \Log = e^* \Log =
    \displaystyle  \prod\limits_{k \geq 0} Sym^k \hh$ (cf. partie \ref{secsc}).

\end{dem}

\subsection{Images directes supérieures du logarithme} \label{secimadir}

\begin{thm}{} \label{imadir} $-$ \\

\begin{itemize}
\item[a)  ]
    On a $H^k \pi_* \Log(d)  =  0$  si $k \not= 2d$.
              Le morphisme $\Log(d) \to \Q(d)$ induit par $\varepsilon \colon \Log \to \Q(0)$
 induit le morphisme
              $$  H^{2d} \pi_* \Log(d) \to
	          H^{2d} \pi_*  \Q(d) = \Q(0).$$
 Ce dernier est un isomorphisme. \\

\item[b) ]  Des deux identités
            $ \displaystyle e^* \Log = \prod\limits_{k \geq 0} Sym^k \hh$
	                (cf. partie \ref{secsc}) et
	       $e^! \Log(d)  =  e^* \Log [-2d]$,

              on déduit que $H^k e^! \Log(d) = 0$, si $k \not= 2d$ et
              $H^{2d} e^! \Log(d)  =  \displaystyle \prod\limits_{k \geq 0} \; \Sym^k \hh $.
              Cette propriété, a) et la suite exacte longue de cohomologie
              associée au triangle distingué:
              $$       e^! \Log(d) \to
	         \pi_* \Log(d) \to  (\pi_U)_* \Log(d)_U \to
	          e^! \Log(d)  [1]           $$
              donnent
              $H^k (\pi_U)_* \Log_U(d) =  0$ si $k \not= 2d-1$ et une suite exacte courte:
	      $$
                     0 \to H^{2d-1} (\pi_{U})_* \Log(d)_U \to
		       H^{2d} e^! \Log(d)  \to  H^{2d} \pi_*
		        \Log(d)  \to 0. $$
              On vérifie que celle-ci s'insère dans le diagramme commutatif suivant:
              $$\xymatrix{ 0 \ar[r] &  H^{2d-1} (\pi_U)_* \Log(d)_U \ar[r]^-{\rho'} &
	             \prod\limits_{k \geq 0} \Sym^k \hh  \ar[r]^-{pr_{k=0}} &
		       \Q(0)  \ar[r] & 0 .\\
		      0 \ar[r] &  H^{2d-1} (\pi_U)_* \Log(d)_U \ar[r] \ar@{=}[u]
		      & H^{2d} e^! \Log(d)  \ar@{=}[u] \ar[r]  & H^{2d} \pi_*
		        \Log(d)
		      \ar[u]^-{a)} \ar[r] & 0}$$
              La factorisation canonique de $\rho'$ à travers
	        $ \displaystyle \prod\limits_{k > 0} \; \Sym^k \hh \hookrightarrow
               \prod\limits_{k \geq 0} \; Sym^k \hh $ donne le morphisme r\'esidu
               $$ \rho :
	         H^{2d-1} (\pi_{U})_* \Log(d)_U \to \prod\limits_{k > 0} \; Sym^k \hh $$
              qui est un isomorphisme.
\end{itemize}
\end{thm}

\begin{dem} Pour a), on renvoie à  [W, I-Cor 4.4], [W, III-Thm 1.3] ou [Ki, Prop 1.1.3].\\

\end{dem}

\section{Le polylogarithme d'un sch\'ema ab\'elien}

\subsection{D\'efinition du polylogarithme d'un sch\'ema ab\'elien} \label{defpol}

Les propri\'et\'es du logarithme \'enonc\'ees dans la partie \ref{proplog}
ont des analogues topologiques \'evidents, e.g.
$R^i (\overline{\pi_U})_* \overline{\Log_U(d)}  = 0$ pour tout $i \not = 2d-1$
et le  morphisme  r\'esidu
$R^{2d-1} (\overline{\pi_U})_* \overline{\Log_U(d)}  \to \overline{e}^* \overline{\Log}$
induit un isomorphisme
$\overline{\rho} \colon R^{2d-1} (\pi_U)_* \overline{\Log_U(d)} \isom
\prod\limits_{n=1}^{\infty} \overline{\Sym^n \hh}$.  \\

On d\'efinit deux isomorphismes $\kappa$ et $\overline{\kappa}$ par le diagramme commutatif, not\'e
$\mathcal{D}_1$, suivant.

$$ \xymatrix{
   & \Ext_{MHM(U)}^{2d-1}( \pi_U^* \hh , \Log_U (d))
   \ar[r]^-{\For} \ar@{=}[d] \ar@{}[dr] |{\text{(adjonction)}}
     \ar `l[ld] `[ddd]_-{\kappa \; }^-{\sim} [ddd]
   &
   \Ext^{2d-1}_{\mathcal{F}_{\Q}(U)}( \overline{\pi_U}^* \overline{\hh} , \overline{\Log_U(d)})
   \ar@{=}[d]
      \ar `r[rd] `[ddd]^-{\; \overline{\kappa}}_-{\sim} [ddd] &
   \\
   & \Ext_{MHM(S)}^{2d-1}(  \hh , (\pi_U)_*\Log_U (d))
   \ar[r]^-{\For} \ar@{=}[d] \ar@{}[dr] |{Thm \ref{imadir}}  &
   \Ext^{2d-1}_{\mathcal{F}_{\Q}(S)}( \overline{\hh} , R \overline{\pi_U}_* \overline{\Log_U(d)})
   \ar@{=}[d]  & \\
   & \Hom_{MHM(S)}( \hh ,    H^{2d-1} (\pi_U)_* \Log_U(d) )
   \ar[r]^-{\For} \ar[d]_-{\rho_*}^-{\sim} \ar@{}[dr] |{\text{(prop. d) de } \Log)} &
   \Hom_{\mathcal{F}_{\Q}(S)}( \overline{\hh} ,   R^{2d-1}  \overline{\pi_U}_* \overline{\Log_U(d)} )
   \ar[d]^-{\overline{\rho}_*}_-{\sim} & \\
   & \Hom_{MHM(S)}( \hh , \prod\limits_{n=1}^{\infty} \Sym^n \hh)
   \ar@{^{(}->}[r]^-{\For} &
   \Hom_{\mathcal{F}_{\Q}(S)}(\overline{\hh},\prod\limits_{n=1}^{\infty}  \Sym^n \overline{\hh}) & \\
}$$
\\

La commutativit\'e du centre de ce diagramme r\'esulte
de la compatibilit\'e du formalisme des 6 foncteurs de $D^b MHM(\cdot)$ et de celui de $D_c^b(\cdot)$
via le foncteur $\For$, e.g. $\For \circ f_* =  R\overline{f}_* \circ \For$ pour $f$ un morphisme
entre sch\'emas de type fini, s\'epar\'es sur $\C$.
On remarque que le but de $\kappa$ s'identifie naturellement \`a $\Hom_{VSHM(S)}( \hh , \hh)$
(cf. pleine fid\'elit\'e de $\iota_S$ et poids).\\

\begin{defi} $-$
Le polylogarithme du sch\'ema ab\'elien $A/S$, not\'e $\Pol$, est d\'efini par
$$   \Ext_{MHM(U)}^{2d-1}( \pi_U^* \hh , \Log_U (d))   \ni \Pol := \kappa^{-1}(Id_{\hh}).$$
\end{defi}

\subsection{Propri\'et\'es du polylogarithme d'un sch\'ema ab\'elien}

\subsubsection{Description compl\`ete dans le cas elliptique ($d=1$)}

Pour tout $\mathbb{V}, \mathbb{W} \in Ob(VSHM(U))$, le foncteur $\iota_U$ induit un
isomorphisme (cf. remarque suivant le Th\'eor\`eme 3.27 de \cite{s})
$$\Ext_{MHM(U)}^{1}( \mathbb{V} ,\mathbb{W} ) \isom \Ext_{VSHM(U)}^{1}( \mathbb{V} ,\mathbb{W} ).$$
Le polylogarithme est une $1$-extension dans $VSHM(U)$ dont une description compl\`ete a \'et\'e donn\'ee
par Beilinson et Levin dans \cite[4.8]{bl}.
On peut \'egalement consulter le th\'eor\`eme \cite[V-Thm 3.4]{w} et sa preuve.

\subsubsection{Sur une description  en dimensions sup\'erieures $(d \geq 2)$}

On d\'emontre que le polylogarithme {\it n'est pas} dans l'image du
morphisme $$\Ext_{VSHM(U)}^{2d-1}( \pi_U^* \hh , \Log_U (d)) \to \Ext_{MHM(U)}^{2d-1}( \pi_U^* \hh , \Log_U (d))$$
induit par $\iota_{S}$ (cf. \cite[III-Thm 2.3 b)]{w}). \\

\subsubsection{Rigidit\'e du polylogarithme d'un sch\'ema ab\'elien}

\begin{lemme} $-$ \label{caracpoltop1}
L'application
$\For \colon \Ext^{2d-1}_{MHM_{\Q}(U)}( \pi_U^* \hh , \Log_U(d)) \to
\Ext^{2d-1}_{\mathcal{F}_{\Q}(U)}( \overline{\pi_U}^* \overline{\hh} , \overline{\Log_U(d)})$
est injective
et $\For(\Pol)$ est caract\'eris\'e par $$\overline{\kappa}(\For(\Pol)) = Id_{\overline{\hh}}.$$
\end{lemme}

\begin{dem} C'est une cons\'equence de la commutativit\'e du  diagramme $\mathcal{D}_1$ et de la
d\'efinition de $\Pol$.

\end{dem}

On a mentionn\'e au d\'ebut de la partie \ref{defpol} que les propri\'et\'es du logarithme
(cf. partie \ref{proplog}) admettent des analogues topologiques. On a alors donn\'e un exemple
en consid\'erant des coefficients rationnels. En fait, ces propri\'et\'es au niveau topologique
peuvent \'egalement se d\'emontrer en consid\'erant des coefficients complexes et on a des r\'esultats
de compatibilit\'es par extension des scalaires de $\Q$ \`a $\C$. Par exemple,
$R^i (\overline{\pi_U})_* \overline{\Log_U}(d)_{\C}  = 0$ pour tout $i \not = 2d-1$
et le  morphisme  r\'esidu
$R^{2d-1} (\overline{\pi_U})_* \overline{\Log_U(d)}_{\C}  \to \overline{e}^* \overline{\Log}_{\C}$
induit un isomorphisme
$\overline{\rho}_{\C} \colon R^{2d-1} (\pi_U)_* \overline{\Log_U(d)}_{\C} \isom
\prod\limits_{n=1}^{\infty} \overline{\Sym^n \hh}_{\C}$.  \\

On d\'efinit un morphisme $\overline{\kappa}_{\C}$ par le diagramme commutatif, not\'e
$\mathcal{D}_2$, suivant

$$ \xymatrix{
   & \Ext^{2d-1}_{\mathcal{F}_{\Q}(U)}( \overline{\pi_U}^* \overline{\hh} , \overline{\Log_U(d)})
   \ar[r]^-{} \ar@{=}[d] \ar@{}[dr] |{\text{(adjonction)}}
     \ar `l[ld] `[ddd]_-{\overline{\kappa} \; }^-{\sim} [ddd]
   &
   \Ext^{2d-1}_{\mathcal{F}_{\C}(U)}( \overline{\pi_U}^* \overline{\hh}_{\C} , \overline{\Log_U(d)}_{\C})
   \ar@{=}[d]
      \ar `r[rd] `[ddd]^-{\; \overline{\kappa}_{\C}}_-{\sim} [ddd] &
   \\
   & \Ext^{2d-1}_{\mathcal{F}_{\Q}(S)}( \overline{\hh} , R \overline{\pi_U}_* \overline{\Log_U(d)})
   \ar[r]^-{} \ar@{=}[d] \ar@{}[dr] |{\text{(prop. du log.)}}  &
   \Ext^{2d-1}_{\mathcal{F}_{\C}(S)}( \overline{\hh}_{\C} , R \overline{\pi_U}_* \overline{\Log_U(d)}_{\C})
   \ar@{=}[d]  & \\
   & \Hom_{\mathcal{F}_{\Q}(S)}( \overline{\hh} ,   R^{2d-1}  \overline{\pi_U}_* \overline{\Log_U(d)} )
   \ar[r]^-{} \ar[d]_-{\rho_*}^-{\sim}  &
   \Hom_{\mathcal{F}_{\C}(S)}( \overline{\hh}_{\C} ,   R^{2d-1}  \overline{\pi_U}_* \overline{\Log_U(d)}_{\C} )
   \ar[d]^-{(\overline{\rho}_{\C})_*}_-{\sim} & \\
   & \Hom_{\mathcal{F}_{\Q}(S)}(\overline{\hh},\prod\limits_{n=1}^{\infty}  \Sym^n \overline{\hh})
   \ar@{^{(}->}[r]^-{} &
   \Hom_{\mathcal{F}_{\C}(S)}(\overline{\hh}_{\C},\prod\limits_{n=1}^{\infty}  \Sym^n \overline{\hh}_{\C}) & \\
}$$

\noindent dans lequel les fl\`eches horizontales sont induites par l'extension des scalaires de $\Q$ \`a $\C$.

\begin{lemme} $-$ \label{caracpoltop2}
L'extension $\For(\Pol) \in
\Ext^{2d-1}_{\mathcal{F}_{\Q}(U)}( \overline{\pi_U}^* \overline{\hh} , \overline{\Log_U(d)})
\subseteq \Ext^{2d-1}_{\mathcal{F}_{\C}(U)}( \overline{\pi_U}^* \overline{\hh}_{\C} , \overline{\Log_U(d)}_{\C})$
est caract\'eris\'ee par $$\overline{\kappa}_{\C}(\For(\Pol)) = Id_{\overline{\hh}_{\C}}.$$
\end{lemme}

\begin{dem} L'inclusion
$\Ext^{2d-1}_{\mathcal{F}_{\Q}(U)}( \overline{\pi_U}^* \overline{\hh} , \overline{\Log_U(d)})
\subseteq \Ext^{2d-1}_{\mathcal{F}_{\C}(U)}( \overline{\pi_U}^* \overline{\hh}_{\C} , \overline{\Log_U(d)}_{\C})$
r\'esulte de la commutativit\'e du diagramme $\mathcal{D}_2$ et de la caract\'erisation du Lemme
\ref{caracpoltop1}.\\

\end{dem}

\subsection{Description du polylogarithme d'un sch\'ema ab\'elien au niveau topologique}

L'objectif de cette partie est de d\'emontrer que les courants d\'efinis par Levin dans \cite{l}
permettent de d\'ecrire
$$\For(\Pol) \in
\Ext^{2d-1}_{\mathcal{F}_{\Q}(U)}( \overline{\pi_U}^* \overline{\hh} , \overline{\Log_U(d)})
\subseteq \Ext^{2d-1}_{\mathcal{F}_{\C}(U)}( \overline{\pi_U}^* \overline{\hh}_{\C} , \overline{\Log_U(d)}_{\C}).$$
On d\'emontre ainsi un r\'esultat qui avait \'et\'e conjectur\'e par Levin.

\subsubsection{\'Equation diff\'erentielle et polylogarithme}

On consid\`ere
le complexe de de Rham des courants sur $A^{\infty}$ \`a valeurs dans le pro-fibr\'e vectoriel
$\G(d)_{\C}$ (cf. partie \ref{courantsplat}),
$(\mathcal{A}^{\bullet} ( \G(d)_{\C}) :=
   \prod\limits_{n=0}^{\infty}
   (\overline{(\Sym^n  \pi^* \hh)(d)}_{\C}) \otimes \mathcal{A}^{\bullet}_{A^{\infty}}
    , \nabla_{\C}^{\bullet}).$
C'est une r\'esolution $\overline{\pi}_*$-acyclique de  $(\overline{\Log(d)})_{\C}$.\\

\begin{nota} $-$ Soit $f \colon \overline{\pi_U}^* \overline{\hh}_{\C} \to
               \mathcal{A}^{2d-1}( \G(d)_{\C})_{|\overline{U}}      $
un morphisme tel que $(\nabla_{\C}^{2d-1})_{|\overline{U}}  \circ f = 0$. 
Le diagramme

$$
\xymatrix{
         &            &   0 \ar[r] & \overline{\pi_U}^* \overline{\hh}_{\C} \ar[r] \ar[d]^-{\; f} & 0  \\
0 \ar[r] & (\G(d)_{\C})_{|\overline{U}} \ar[r]^-{(\nabla_{\C})_{|\overline{U}} } & \ldots \ar[r]
         &  \mathcal{A}^{2d-1} (\G(d)_{\C})_{|\overline{U}} \; \ar[r]^-{(\nabla_{\C}^{2d-1})_{|\overline{U}}} &
           \;  \mathcal{A}^{2d} (\G(d)_{\C})_{|\overline{U}} \ar[r] & \ldots \\
0 \ar[r] & (\overline{\Log_U (d)})_{\C} \ar[u]^-{qis \;} \ar[r] & 0 \\
         }$$
d\'efinit un \'el\'ement de
$\Hom_{D^b(\mathcal{F}_{\C}(U))}(\overline{\pi_U}^* \overline{\hh}_{\C},
                       (\overline{\Log_U(d)})_{\C}[2d-1])$
que l'on note $M(f)$.
\end{nota}

\begin{thm} $-$ \label{pol_equa_diff}
Soit $f \colon \overline{ \pi^* \hh }_{\C} \to \aaa^{2d-1}_{A^{\infty}}( \G(d)_{\C})$ un morphisme dans
 $\F_{\C}(A)$ v\'erifiant la propri\'et\'e $(P)$ suivante:
 $$ (P) \quad \quad \quad \nabla^{2d-1} \circ f = 
      (2 \pi i)^d \; \delta_{S^{\infty}} \; Id_{\overline{ \pi^* \hh }_{\C}},$$
 o\`u $S^{\infty}$ est vue comme une sous-vari\'et\'e ferm\'ee de $A^{\infty}$ via $e^{\infty}$.     
 Alors, on a:
 \begin{center}
 \begin{tabular}{ll}
   1. & $ \nabla^{2d-1}_{|\overline{U}} \circ f_{|\overline{U}} = 0$. \\
   2. &  $M(f_{|\overline{U}}) = \overline{\Pol}.$
 \end{tabular}
 \end{center}        
\end{thm}

\begin{notas} Pour $(E,\nabla)$ un fibr\'e vectoriel r\'eel \`a connexion plate et
$f \colon (E_1,\nabla_1) \to  (E_2,\nabla_2)$ un morphisme de
fibr\'es vectoriels r\'eels \`a connexions plates, on note: \\

\begin{tabular}{cp{15cm}}
$E^{\circ}$ & le syst\`eme local $\Ker(\nabla)$,\\
$f^{\circ}$ & le morphisme de syst\`emes locaux induit par $f$ entre $E_1^{\circ}$ et $E_2^{\circ}$.
\end{tabular}

\end{notas}

\begin{dem} \\

\begin{itemize}
 \item[1. ] C'est une cons\'equence imm\'ediate de la propri\'et\'e $(P)$. \\

 \item[2. ] D'apr\`es le Lemme \ref{caracpoltop2}, il suffit de d\'emontrer l'assertion suivante

$$ (A_1) \quad \quad \quad \overline{\kappa}_{\C} ( M(f_{|\overline{U}})  ) = Id_{\overline{\hh }_{\C}},$$
 ce que l'on fait ci-dessous. \\

\begin{itemize}
\item[a)]
On commence par r\'eduire le calcul  de $\overline{\kappa}_{\C} ( M(f_{|\overline{U}})  )$
 dans lequel interviennent des courants \`a valeurs dans 
un pro-fibr\'e vectoriel $\G(d)_{\C}$ \`a plusieurs calculs ne mettant en jeu que des courants \`a valeurs 
 dans des fibr\'es vectoriels (les fibr\'es $\G_l(d)_{\C}$, tronqu\'es de $\G(d)_{\C}$). On rappelle que le morphisme $\overline{\kappa}_{\C}$ est donn\'e par la composition
$$ (\overline{\rho}_{\C})_* \circ  H^0 \circ  adj, $$
o\`u \\

\begin{tabular}{p{15cm}}
 
$adj \colon \Hom_{D^b \F_{\C}(U)} ( \overline{\pi_U^* \hh}_{\C} , \overline{\Log_U(d)}_{\C}[2d-1]) \to
         \Hom_{D^b \F_{\C}(S)} ( \overline{\hh}_{\C} , R (\overline{\pi_U})_*   
	                       \overline{\Log_U(d)}_{\C}[2d-1] )$ est l'isomorphisme d'adjonction,   \\
	                                              
			       \\
 $H^0 \colon \Hom_{D^b \F_{\C}(S)} ( \overline{\hh}_{\C} , R (\overline{\pi_U})_*   
	                       \overline{\Log_U(d)}_{\C}[2d-1] ) \to
        \Hom_{\F_{\C}(S)} ( \overline{\hh}_{\C} , R^{2d-1}(\overline{\pi_U})_*   
	                       \overline{\Log_U(d)}_{\C})$,  \\
	                       
	                       	\end{tabular}                       
	                       
\begin{tabular}{p{15cm}}	
			       \\
 $\overline{\rho}_{\C} \colon    R^{2d-1}(\overline{\pi_U})_* \overline{\Log_U(d)}_{\C} \to
                                \overline{e}^* \overline{\Log}_{\C}$
est un morphisme de bord qui appara\^it dans la suite exacte longue de cohomologie locale associ\'ee 
\`a la situation g\'eom\'etrique suivante:
$$  \xymatrix{
     \overline{e} \colon \overline{S} \; \ar@{^{(}->}[r]|-{|} & \overline{A}  &  \; \overline{U}  
                              \ar@{_{(}->}[l]|-{\circ}},$$
avec comme coefficient le pro-syst\`eme local $\overline{\Log(d)}_{\C}$.
  			       \\
			       \\
\end{tabular}
Pour prouver la relation  
$\overline{\kappa}_{\C} ( M(f_{|\overline{U}})  ) = Id_{\overline{\hh }_{\C}}$, 
il suffit de d\'emontrer que pour tout $l \in \N^{\geq 2}$, l'assertion suivante est valide:
$$ (A_2^l) \quad \quad \quad
\overline{\kappa}_{\C}^l ( M^l ( ( p_l^c  \circ f)_{\overline{U}})) = Id_{\overline{\hh }_{\C}},$$

o\`u \\ 

\begin{tabular}{p{15cm}}
$\overline{\kappa}_{\C}^l \colon  
 \Hom_{D^b \F_{\C}(U)} ( \overline{\pi_U^* \hh}_{\C} , (\G_l^{\circ})_{\overline{U}} (d)_{\C} [2d-1])
 \to \overline{e}^* (\G_l^{\circ})_{\C} = \prod\limits_{k =1}^l \Sym^k \overline{\hh }_{\C}$
 est d\'efini de mani\`ere analogue \`a $\overline{\kappa}_{\C}$ en prenant cette fois
 $\G_l^{\circ}(d)_{\C}$ comme coefficient, \\
 \\
$M^l( ( p_l^c  \circ f)_{\overline{U}})$ est d\'efini de mani\`ere analogue \`a 
       $M(f_{|\overline{U}})$ en consid\'erant  $\G_l^{\circ}(d)_{\C}$ comme coefficient 
       (cf. ci-dessous),\\ \\
 
$ p^c_l \colon \aaa^{2d-1}_{A^{\infty}} (\G(d)_{\C}) \to 
                            \aaa^{2d-1}_{A^{\infty}} ( \G_l (d)_{\C} )$ est
 le morphisme obtenu en poussant les $(2d-1)$-courants sur $A^{\infty}$ \`a valeurs dans 
 $\G(d)_{\C}$ \`a l'aide du morphisme $p_l$ dans la partie 3.3.2.
 \\
 \\
\end{tabular}

On fixe $l \in \N^{\geq 2}$ pour la suite de la d\'emonstration. Puisque $f$ v\'erifie la propri\'et\'e
$(P)$, $ p_l^c \circ f $ v\'erifie la propri\'et\'e suivante:
$$  \nabla_l^{2d-1} \circ p_l^c \circ f  =
      (2 \pi i)^d \; \delta_{S^{\infty}} \; Id_{\overline{ \pi^* \hh }_{\C}}$$
et en particulier $(\nabla_l^{2d-1} \circ p_l^c \circ f)_{|\overline{U}} = 0$, 
ce qui implique que   $M^l( ( p_l^c  \circ f)_{\overline{U}})$ est bien d\'efini.  \\

\item[b)] 
L'assertion $(A^l_2)$ est de nature locale. Soient $s \in S(\C)$ et
$V$ un voisinage ouvert connexe et simplement connexe de $s$ dans $S(\C)$. On souhaite d\'ecrire
le morphisme 
$$ \overline{\kappa}_{\C}^l ( M^l ( ( p_l^c  \circ f)_{\overline{U}}))_{V}
   \colon \Gamma(V, \overline{\hh}_{\C}) \to  \Gamma(V, \overline{\hh}_{\C}).$$

On introduit, pour ce faire, la notation suivante.
Soient $X$ une vari\'et\'e alg\'ebrique complexe, 
       $F$ une sous-vari\'et\'e alg\'ebrique ferm\'ee de codimension pure $d$,
       $i \colon F \hookrightarrow X$ l'immersion ferm\'ee correspondante
       et $\mathbb{V}$ un (pro-)syst\`eme local de $\R$-vectoriels sur $\overline{X}$.
Alors on a une identification canonique 
$\overline{i}^{!} \mathbb{V} = \overline{i}^*V (-d) [-2d]$.
Soient $F'$ un ouvert de $\overline{F}$ et $X'$ un ouvert de $\overline{X}$ contenant $F'$.  
On note 
$$ \rho(F',X',\mathbb{V}) \colon H^{2d-1}(X' \setminus F' , \mathbb{V}(d)_{\C}) \to  
                        H^{2d}( F' , \overline{i}^! \mathbb{V}(d)_{\C} )   = 
			 \Gamma(F', \overline{i}^* \mathbb{\V}_{\C})$$
le morphisme de bord qui appara\^it dans la suite exacte longue de cohomologie locale.\\

Le  $\overline{\kappa}_{\C}^l ( M^l ( ( p_l^c  \circ f)_{\overline{U}}))_{V}$
 est donn\'e 
par la composition suivante:

$$\xymatrix{  \Gamma(V, \overline{\hh}_{\C}) \ar@{=}[d] \\
             \Gamma(\overline{\pi}^{-1}(V), \overline{\pi^*\hh}_{\C}) 
	    \ar[d]^-{ \; (p_l^c \circ f)_{\overline{\pi}^{-1}(V)}} \\
\left\{ 
c \in \Gamma( \overline{\pi}^{-1}(V) , \aaa^{2d-1}_{A^{\infty}}( (\G_l(d))_{\C}) )
      \; : \;  (\nabla_l^{2d-1} (c))_{| \overline{\pi}^{-1}(V) \setminus V }  = 0 \right\}  
\ar[d] \\
H^{2d-1}( \overline{\pi}^{-1}(V) \setminus V ,  (\G^{\circ}_l(d))_{\C}) ) 
\ar[d]^-{  \;  \rho( V , \overline{\pi}^{-1}(V), \G^{\circ}_l)} \\
\Gamma( V , \overline{e}^*  (\G^{\circ}_l)_{\C}) \ar@{=}[d] \\
\prod\limits_{k=0}^l  \Sym^k \Gamma( V ,  \overline{\hh}_{\C} ).}$$  %
 
Pour d\'emontrer localement en $s$ l'assertion $(A^l_2)$, il suffit donc de prouver que: \\

\begin{center} 
\begin{tabular}{cp{11cm}}
$(A^l_3)$ &
 pour tout
$ h \in   \Gamma(V, \overline{\hh}_{\C})$, 
$c \in \Gamma( \overline{\pi}^{-1}(V) , \aaa^{2d-1}_{A^{\infty}}( \G_l(d)_{\C}) )$
tel que  $\nabla_l^{2d-1} (c) = (2 \pi i )^d \; \delta_V \; h$, on a
$ \rho( V , \overline{\pi}^{-1}(V), \G^{\circ}_l ) ([c_{|\overline{\pi}^{-1}(V) \setminus V}]) = h.$\\
\\
\end{tabular}
\end{center}

\item[c)] On explique maintenant comment passer du coefficient $\G^{\circ}_l$ au coefficient
          trivial $\R$.
          Le morphisme  $ \rho( V , \overline{\pi}^{-1}(V), \G^{\circ}_l )$ \'etant un 
          morphisme de bord dans une suite exacte de cohomologie locale, on peut remplacer
	  $\overline{\pi}^{-1}(V)$ par un voisinage ouvert de $V$ dans $\overline{A}$. 
	  Soit $W$ un voisinage ouvert de $e(s)$ dans $\overline{A}$ sur lequel le fibr\'e
	  vectoriel r\'eel \`a connexion int\'egrable $(\G_l,\nabla_l)$ est isomorphe 
	  au fibr\'e trivial de fibre
	  $\prod\limits_{k=0}^l  \Sym^k \Gamma( V ,  \overline{\hh}_{\R} )$
	  muni de la connexion de Gauss-Manin. Quitte \`a remplacer $V$ 
	  par un voisinage ouvert de $s$ dans $\overline{S}$ qui est connexe et simplement connexe, 
	  on peut supposer que $V  \subseteq W$. On est ainsi ramen\'e au cas 
	  o\`u le  coefficient est  trivial, i.e. il suffit de d\'emontrer que:\\

\begin{center}
\begin{tabular}{cp{11cm}}
$(A_4)$ &
 pour tout $c \in \Gamma( W , \aaa^{2d-1}_{A^{\infty}}( \R(d)_{\C}) )$
tel que  $d c = (2 \pi i )^d \; \delta_V $, on a
$ \rho( V , W , \R) ([c_{|W \setminus V}]) = 1.$\\
\\
\end{tabular}
\end{center}

\item[d)] Compte-tenu du caractère local de l'assertion $(A_4)$ et
de la structure locale des immersions ferm\'ees en g\'eom\'etrie analytique,
il suffit, modulo l'application d'un biholomorphisme, de démontrer l'assertion $(A_4)$ dans la situation g\'eom\'etrique suivante: \\

\begin{center}
\begin{tabular}{ll}
i) & $V$ est une boule ouverte de $\C^{n}$ contenant 0 ($n \in \N$). \\
ii) & $ s = 0 \in V$. \\
iii) & $W = V \times B(0,1)$ o\`u $B(0,1)$ est la boule ouverte de $\C^{d}$ centr\'ee en $0$ et de rayon 1. \\
iv) & L'immersion ferm\'ee
$e^{\infty} \colon V \hookrightarrow V \times B(0,1)$ est donn\'ee par $v \mapsto (v,0)$.\\
& \\
\end{tabular} 
\end{center}

On se place d\'esormais dans ce contexte g\'eom\'etrique. 
On souhaite maintenant r\'eduire la d\'emonstration de l'assertion
$(A_4)$ \`a la preuve d'un cas particulier de celle-ci: $V=\{0\}$,
$W=B(0,1)$ et $c$ est le courant associ\'e \`a la forme de Bochner-Martinelli
dont on rappelle succinctement la construction. On donne \'egalement 
l'\'equation diff\'erentielle que satisfait ce courant.  \\

Soit $\beta$ la $(2d-1)$-forme diff\'erentielle sur $B(0,1) \setminus \{0\}$
$$ \beta := F^* K,$$
o\`u $K$ d\'esigne le noyau de Bochner-Martinelli et $F$ est l'application de
$B(0,1) \setminus \{0\}$ dans $\C^d \times \C^d$ d\'efinie par
$F(z) = (2z,z)$ pour $z \in B(0,1) \setminus \{0\}$ (cf. \cite[p. 371 et 655]{gh}).
	  Alors $d\beta = 0$ et
	  les coefficients de $\beta$ sont localement $L^1$. Ainsi,
	  $\beta$ d\'efinit un courant sur $B(0,1)$ que l'on note $\overline{\beta}$.
	  La d\'eriv\'ee de ce courant v\'erifie
	  $d \overline{\beta} = \delta_0$  (cf. \cite[p. 371 et 372]{gh}). \\

On introduit alors l'assertion suivante  
\begin{center}
\begin{tabular}{cp{11cm}}
& \\
$(A_5)$ &
$ \rho( \{ 0 \} , B(0,1) , \R) ((2 \pi i)^d \beta) = 1.$\\
\\
\end{tabular}
\end{center}
et on d\'emontre que celle-ci implique $(A_4)$.\\

On suppose l'assertion $(A_5)$ v\'erifi\'ee et on fixe
$c \in \Gamma( V \times B  , \aaa^{2d-1}_{A^{\infty}}( \R(d)_{\C}) )$
tel que  $d c = (2 \pi i )^d \; \delta_V $. \\

\begin{itemize}
\item[$\bullet$] Soit $pr \colon V \times B(0,1) \to B(0,1)$ la projection canonique. 
On remarque, en consid\'erant l'expression en coordonn\'ees de $\beta$ que tous
les coefficients de $(pr^{|B(0,1) \setminus \{0\}})^* \beta$ sont $L^1$. 
La forme diff\'erentielle $(pr^{|B(0,1) \setminus \{0\}})^* \beta$ d\'efinit donc 
un courant que l'on note $pr^* \overline{\beta}$. (L'existence d'un pullback n'est pas assur\'ee pour les courants en gen\'eral et c'est cette propri\'et\'e d'extension de la forme 
$(pr^{|B(0,1) \setminus \{0\}})^* \beta$ en un courant d\'efini sur $W$
qui explique, entre autre, la consid\'eration de $\beta$. Une autre motivation est la formule de Bochner-Martinelli utilis\'ee ci-apr\`es.) 
D'autre part, on v\'erifie, \`a l'aide de l'\'equation  
 $d \overline{\beta} = \delta_0$ que le courant $pr^* \overline{\beta}$
	satisfait l'équation:
	 $$ (*) \quad \quad \quad d \; pr^* \overline{\beta} = \delta_V.$$

\item[$\bullet$] Comme $\rho( V , W , \R)$ est le morphisme bord
	  d'une suite exacte longue de cohomologie locale et
	  $dc = (2 \pi i)^d d \; pr^* \overline{\beta}$ (d'apr\`es (*)), on a

$$ \begin{array}{llll}  &
	    \rho( V , W , \R) ([c_{|W \setminus V}]) & = &
             \rho( V , W , \R) ([ (2 \pi i)^d  (pr^* \overline{\beta})_{|W \setminus V}]) \\
 & & = & (2 \pi i)^d \rho( V , W , \R) ( [   pr^{|B(0,1) \setminus \{0\}})^* \beta]).
\end{array}$$

          Il suffit donc de considérer le cas particulier $c =(2 \pi i)^d  pr^* \overline{\beta}$
	  pour démontrer $(A_4)$. \\

\item[$\bullet$] On consid\`ere le diagramme suivant:

$$ \xymatrix{ V \ar@{^{(}->}[r]^-{e^{\infty}}
             & V \times B(0,1)     \\
             \{ 0 \}  \ar@{^{(}->}[r]_-{(e^{\infty})_0} \ar@{^{(}->}[u]^-{i}
	     &
             B(0,1)  \ar@{^{(}->}[u]_-{i'} }$$
dans lequel les morphismes $i$, $i'$ et $(e^{\infty})_0$ sont d\'efinis par
$$i(0)=0, \quad \; \forall \; b \in B(0,1) \; \; i'(b)=(0,b)
\quad \text{ et } \quad (e^{\infty})_0(0)=0.$$

On v\'erifie que l'on a la relation suivante:

$$ \begin{array}{llll}  &
 \rho( V , W , \R) ( [(pr_{|B(0,1) \setminus \{0\}})^* \beta ])
   & = &
  \rho( \{0 \} , B(0,1) , \R ) ( ( i'_{|B(0,1) \setminus \{ 0 \} })^*  [ (pr^{|B(0,1) \setminus \{0\}})^* \beta]) \\
 & & = & \rho( \{0 \} , B(0,1) , \R ) ( [ ( i'_{|B(0,1) \setminus \{ 0 \} })^*   (pr^{|B(0,1) \setminus \{0\}})^* \beta]) \\
  & & = & \rho( \{0 \} , B(0,1) , \R )([\beta]). \\
  \end{array}$$

  La preuve de l'implication $(A_5) \imp (A_4)$ est ainsi achevée. \\
\end{itemize}

\item[e)] Il reste donc \`a d\'emontrer l'assertion $(A_5)$.
	  D'apr\`es \cite[V.7]{i}, on a:

	  $$
	    \rho( \{ 0 \} , B(0,1) , \R) ([\beta]) = (2 \pi i)^{-d} \;
	   \int_{\partial B(0,r)} \; \beta_{|\partial B(0,r)},$$
	  o\`u $B(0,r)$ est la boule de $\C^d$ centr\'ee en $0$ et de rayon $r \in \; ]0,1[$.
	  L'assertion $(A_5)$ est alors cons\'equence de la formule de Bochner-Martinelli
           (cf. \cite[p. 372]{gh}):

$$ \int_{\partial B(0,r)} \; \beta_{|\partial B(0,r)} =  1. $$

 \end{itemize}

  \end{itemize}

\end{dem}

\subsubsection{Les courants de Levin} \label{courantslevin}

Soit $\omega$ une polarisation
du sch\'ema ab\'elien $\pi \colon A \to S$.
Dans \cite{l},  Levin d\'efinit, \`a partir de $\omega$, des s\'eries de formes diff\'erentielles sur $A^{\infty}$ \`a valeurs dans $\Sym^{a-1} (\OO_{A^{\infty}} \otimes \overline{\pi^* \hh_{\C}})$, not\'ees $g'_a$ $(a \in \N^*)$. On pr\'ecise ci-dessous en quel sens ces s\'eries convergent et on donne des indications quant \`a la mani\`ere d'\'etablir ces r\'esultats de convergence.\\

\begin{itemize}
\item[$\bullet$] 
Pour $a > 2d$, $g'_a$ converge uniform\'ement vers une forme diff\'erentielle
sur $A^{\infty}$. \\

En effet, comme il s'agit d'une assertion de nature locale, on peut supposer que $S$ est un ouvert de $\C^n$ et que $(\pi \colon A \to S, \omega)$ est un pullback de la famille universelle de vari\'et\'es ab\'eliennes polaris\'ee consid\'er\'ee par Levin
(cf. \cite[2.3]{l}), en modifiant \'eventuellement la polarisation qu'il introduit, de mani\`ere \`a tenir compte du type de la polarisation $\omega$. Dans ce cas, on dispose de coordonn\'ees globales et d'une formule "explicite" pour $g'_a$. On montre alors la convergence de $g'_a$ en utilisant que la s\'erie num\'erique
$$ \sum_{n_1 + \dots + n_{2d}  \in \Z^{2d} \setminus \{0\}} \; \;
   (n_1^2 + \dots + n^2_{2d})^{-a/2} $$
converge si $a > 2d$.  \\ 

\item[$\bullet$] 
 Pour $a \leq 2d$, $g'_a$, vues comme s\'eries de courants, convergent au sens des courants (cf. partie \ref{convergencecourants}). \\

Pour le voir, l'\'enonc\'e \'etant local, on peut proc\'eder comme ci-dessus pour obtenir une formule "explicite" de $g'_a$. On applique alors 
un op\'erateur de Laplace (associ\'e aux coordonn\'ees verticales), \'eventuellement plusieurs fois,  \`a des s\'eries de formes diff\'erentielles convergeant uniform\'ement (dont la convergence peut s'\'etablir comme celle des s\'eries $g'_a$ pour $a > 2d$) pour obtenir $g'_a$
et conclure.
La d\'emonstration est analogue \`a celle de \cite[Thm 3 1.3]{t}. \\ 
\end{itemize}

\`A l'aide de ces s\'eries, il construit un morphisme
$\PP_{\omega} \colon \overline{\pi}^* \overline{\hh}_{\C} \to \displaystyle
 \mathcal{A}^{2d-1} (\G(d)_{\C})$
(cf. \cite[Thm 3.4.4]{l}).
Pour une expression explicite de $\PP_{\omega}$, dans le cas o\`u le sch\'ema ab\'elien est une
famille modulaire de Siegel (resp. Hilbert-Blumenthal), on peut consulter \cite[2.3]{l}
(resp. \cite{b}). Levin d\'emontre que $\PP_{\omega}$ v\'erifie la propri\'et\'e $(P)$ du th\'eor\`eme pr\'ec\'edent
\cite[Thm 3.4.4]{l} et conjecture que ce morphisme d\'ecrit $\Pol$. Du Th\'eor\`eme \ref{pol_equa_diff},
on d\'eduit une preuve de cette conjecture. Pr\'ecis\'ement, on a le corollaire suivant.\\

\begin{coro} $-$ \label{explicitation_Pol}
Soit $\omega$ une polarisation du sch\'ema ab\'elien $A/S$.
 Le morphisme $\PP_{\omega}$ de Levin d\'ecrit le polylogarithme
au niveau topologique, i.e. $M((\PP_{\omega})_{|\overline{U}}) =  For(\Pol)$.
 $For(\Pol)$ co\"incide donc avec l'\'el\'ement de  
$\Hom_{D^b(\mathcal{F}_{\C}(U))}(\overline{\pi_U}^* \overline{\hh}_{\C} , (\overline{\Log_U(d)})_{\C}[2d-1])$ d\'efini par le diagramme suivant:

$$
\xymatrix{          &            &   0 \ar[r] & \overline{\pi_U}^* \overline{\hh}_{\C} \ar[r] \ar[d]^-{\;(\PP_{\omega})_{|\overline{U}} } & 0  \\
0 \ar[r] & (\G(d)_{\C})_{|\overline{U}} \ar[r]^-{(\nabla_{\C})_{|\overline{U}} } & \ldots \ar[r]
         &  \mathcal{A}^{2d-1} (\G(d)_{\C})_{|\overline{U}} \; \ar[r]^-{(\nabla_{\C}^{2d-1})_{|\overline{U}}} &
           \;  \mathcal{A}^{2d} (\G(d)_{\C})_{|\overline{U}} \ar[r] & \ldots \\
0 \ar[r] & (\overline{\Log_U (d)})_{\C} \ar[u]^-{qis \;} \ar[r] & 0. }
$$
\\
\end{coro}

  On termine cette partie avec un r\'esultat concernant la lissit\'e des courants de 
Levin.

\begin{prop} \emph{(Levin)}  $-$ 
Pour tout ouvert $V$ de $\overline{U}$, 
tout $h \in \Gamma(V,\overline{\pi}^* \overline{\hh}_{\C})$, le courant
le courant $ \PP_{\omega}(h)$ est lisse sur $V$. 
\end{prop}

\begin{dem} L'assertion se déduit de \cite[Proposition 3.4.2]{l} et de la Proposition A2.1 de l'appendice.
\end{dem}

\section{Les classes d'Eisenstein d'un sch\'ema ab\'elien}

Soit $x \colon S \to A$ une section de torsion et soit $l \in \N$. On d\'efinit deux applications
$val_x^l$ et $\overline{val_x^l}$ par le diagramme commutatif (cf. compatibilit\'e des formalismes des 
6 foncteurs au niveau des modules de Hodge et au niveau topologique via le foncteur $For$) suivant not\'e $\mathcal{D}_3$.

$$\xymatrix{
   &  \Ext^{2d-1}_{MHM(U)}( \pi_U^* \hh , \Log_U(d))
     \ar `l[ld] `[dddd]_-{val_x^l} [dddd]
    \ar[r]^-{\For} \ar[d]_-{x^*} \ar@{}[dr] |{\text{(cf. partie \ref{principescindage})}} &
    \Ext^{2d-1}_{\mathcal{F}_{\Q}(U)}( \overline{\pi_U^* \hh} , \overline{\Log_U(d)})
     \ar[d]^-{\overline{x}^*}
    \ar `r[rd] `[dddd]^-{\overline{val_x^l}} [dddd]
    & \\
   &  \Ext^{2d-1}_{MHM(S)}( \hh , \prod\limits_{n=0}^{\infty} (\Sym^n \hh)(d))
    \ar[r]^-{\For} \ar@{}[dr] |{\text{(dualit\'e)}} \ar@{=}[d] &
    \Ext^{2d-1}_{\mathcal{F}_{\Q}(S)}( \overline{\hh} , \prod\limits_{n=0}^{\infty} \overline{(\Sym^n \hh)(d)} )
    \ar@{=}[d] & \\
    & \Ext^{2d-1}_{MHM(S)}( \Q(0) , \prod\limits_{n=0}^{\infty} (\Sym^n \hh) \otimes \hh^{\vee}(d))
    \ar[r]^-{\For} \ar[d]_-{pr_{l+1}} &
    H^{2d-1}_{\text{Betti}} ( \overline{S} ,
     \prod\limits_{n=0}^{\infty} \overline{(\Sym^n \hh)} \otimes \overline{\hh^{\vee}(d))}  )
    \ar[d]^-{pr_{l+1}} & \\
    & \Ext^{2d-1}_{MHM(S)}( \Q(0) , (\Sym^{l+1} \hh) \otimes \hh^{\vee}(d))
    \ar[d]_-{\; \text{contraction}} \ar[r]^-{\For} &
    H^{2d-1}_{\text{Betti}} ( \overline{S} ,
    \overline{(\Sym^{l+1} \hh)} \otimes \overline{\hh^{\vee}(d))})
     \ar[d]^-{\; \text{contraction}} & \\
    & \Ext^{2d-1}_{MHM(S)}( \Q(0) , (\Sym^l \hh) (d))
    \ar[r]^-{\For} &
     H^{2d-1}_{\text{Betti}} ( \overline{S} , \overline{(\Sym^l \hh) (d)} )  & }$$
\\

\begin{defi} $-$
L'extension $val_x^l(\Pol)$ est appel\'ee
$l$-i\`eme classe d'Eisenstein du sch\'ema ab\'elien $\pi \colon A \to S$
associ\'ee \`a $x$ et not\'ee  $\Eis_x^l$.
\end{defi}

\begin{rema} $-$
D'apr\`es un th\'eor\`eme de Kings, $\Eis_x^l$ est d'origine motivique (voir \cite{k}).
\end{rema}

L'application $\overline{val_x^l}$ a un analogue pour des coefficients complexes que l'on d\'efinit
par le diagramme commutatif suivant not\'e $\mathcal{D}_4$ .

$$\xymatrix{
   & \Ext^{2d-1}_{\mathcal{F}_{\Q}(U)}( \overline{\pi_U^* \hh} , \overline{\Log_U(d)})
     \ar `l[ld] `[dddd]_-{\overline{val_x^l}} [dddd]
    \ar[r]^-{ } \ar[d]_-{x^*} \ar@{}[dr] |{\text{(cf. partie \ref{principescindage})}} &
    \Ext^{2d-1}_{\mathcal{F}_{\C}(U)}( \overline{\pi_U^* \hh}_{\C} , \overline{\Log_U(d)}_{\C})
     \ar[d]^-{\overline{x}^*}
    \ar `r[rd] `[dddd]^-{(\overline{val_x^l})_{\C}} [dddd]
    & \\
   & \Ext^{2d-1}_{\mathcal{F}_{\Q}(S)}( \overline{\hh} , \prod\limits_{n=0}^{\infty} \overline{(\Sym^n \hh)(d)} )
    \ar[r]^-{ } \ar@{}[dr] |{\text{(dualit\'e)}} \ar@{=}[d] &
    \Ext^{2d-1}_{\mathcal{F}_{\C}(S)}( \overline{\hh}_{\C} , \prod\limits_{n=0}^{\infty}
    \overline{(\Sym^n \hh)(d)}_{\C} )
    \ar@{=}[d] & \\
    & H^{2d-1}_{\text{Betti}} ( \overline{S} ,
     \prod\limits_{n=0}^{\infty} \overline{(\Sym^n \hh)} \otimes \overline{\hh^{\vee}(d))}  )
    \ar[r]^-{ } \ar[d]_-{pr_{l+1}} &
    H^{2d-1}_{\text{Betti}} ( \overline{S} ,
     \prod\limits_{n=0}^{\infty} \overline{(\Sym^n \hh)}_{\C} \otimes \overline{\hh^{\vee}(d))}_{\C}  )
    \ar[d]^-{pr_{l+1}} & \\
    & H^{2d-1}_{\text{Betti}} ( \overline{S} ,
    \overline{(\Sym^{l+1} \hh)} \otimes \overline{\hh^{\vee}(d))})
    \ar[d]_-{\; \text{contraction}} \ar[r]^-{ } &
    H^{2d-1}_{\text{Betti}} ( \overline{S} ,
    \overline{(\Sym^{l+1} \hh)}_{\C} \otimes \overline{\hh^{\vee}(d))}_{\C})
     \ar[d]^-{\; \text{contraction}} & \\
    & H^{2d-1}_{\text{Betti}} ( \overline{S} , \overline{(\Sym^l \hh) (d)} )
    \ar@{^{(}->}[r]^-{ } &
     H^{2d-1}_{\text{Betti}} ( \overline{S} , \overline{(\Sym^l \hh) (d)}_{\C} )  & }$$

\noindent dans lequel les fl\`eches horizontales sont induites par l'extension des scalaires de
$\Q$ \`a $\C$. \\

\begin{rema} $-$
\'Etant donn\'ee une polarisation $\omega$ du sch\'ema ab\'elien $A/S$, on peut alors expliciter
 $$ For(\Eis_x^l)  \in H^{2d-1}_{\text{Betti}} ( \overline{S} , \overline{(\Sym^l \hh) (d)} )
    \subseteq H^{2d-1}_{\text{Betti}} ( \overline{S} , \overline{(\Sym^l \hh) (d)}_{\C} ) $$
\`a l'aide de l'identit\'e
$$ For(\Eis_x^l)=(\overline{val_x^l})_{\C} (M((\PP_{\omega})_{|\overline{U}}) )$$
qui se d\'eduit du Corollaire \ref{explicitation_Pol} et de la commutativit\'e des diagrammes $\mathcal{D}_3$ et $\mathcal{D}_4$, dans le cas o\`u $l > 2d$ (cf. convergence des s\'eries de Levin discut\'ee dans la partie \ref{courantslevin}). \end{rema}

Dans \cite{b}, on effectue ce calcul pour le sch\'ema ab\'elien universel au-dessus
d'une vari\'et\'e de Hilbert-Blumenthal. Le r\'esultat est que, dans ce cas, $For(\Eis_x^l)$
s'exprime \`a l'aide de s\'eries d'Eisenstein-Kronecker et l'on d\'emontre, en utilisant le Corollaire \ref{explicitation_Pol}, que certaines classes d'Eisenstein sont non nulles en \'etablissant qu'elles d\'eg\'en\`erent au bord de la compactification de Baily-Borel de la base en des valeurs sp\'eciales de fonctions $L$ associ\'ees
au corps de nombres totalement r\'eel sous-jacent.

\newpage

\begin{center}
{\huge  {\bf Appendix}}
 \\ 
 \vspace{0.3cm}
 {\large by} \\
\vspace{0.3cm}
{\Large Andrey Levin} \\
 \vspace{0.5cm}
email: \verb+alevin@wave.sio.rssi.ru+\\
$\;$\\
\vspace{0.8cm}
\end{center}

We prove a smoothness result for the polylogarithmic current defined in \cite{l}. 
Essentially this is rather standard exercice in the Riemann method for analytic continuation of the $\zeta$-function. \\

\noindent {\bf Notations} $-$ For $X$ a complex analytic variety, we denote \\

\begin{tabular}{cp{14cm}}
 $X^{\infty}$ & the $\cinf$ differential variety associated to $X$, \\
 $\OO_{X^{\infty}}$ & the sheaf of real valued differentiable functions on $X^{\infty}$,\\
 ${T}X^{\infty}$ & the real tangent bundle of $X^{\infty}.$
 \end{tabular}

\section*{A1. Polylogarithmic currents}

Let $S$ be a complex analytic variety and $( \pi \colon X \to S, e \colon S \to B , \omega)$ be a family  
of abelian varieties over $S$ as defined in  \cite[1.1.2]{l}, i.e. 
 $\pi$ is a proper smooth morphism of complex analytic varieties of relative dimension $d$, 
 $X_s := \pi^{-1}(s)$ is a $d$-dimensional complex torus for each $s \in S$, 
$e$ is a section of $\pi$ and
$\omega$ is a $(1,1)$-cohomology class on $X^{\infty}$ such that the restriction to each $X_s$
                         is a polarization for each $s \in B$. \\                  

Let $\Lambda$ be the dual of the local system $\R^1 \pi_* \Z$ over $S$. 
Its stalk at $s \in S$ is $H_1(X_b,\Z)$ and it is equipped with a natural structure of variation of pure Hodge structures
of type $\{(-1,0),(0,-1)\}$ over $S$. Thus the complex  vector bundle $\hh:=\Lambda \otimes \OO_{S^{\infty}} \otimes \C$ 
has a canonical Hodge decomposition
$ \hh = \hh^{-1,0} \oplus \hh^{0,-1}$. \\

The polylogarithmic current is a $(2d-2)$-current on $X^{\infty}$ with values in the complex vector bundle 
$ \displaystyle \prod_{k \geq 0} \Sym^k \pi^*\hh$. 
In this part we recall the definition of these currents when $S$ is simply connected. 
For arbitrary $S$, the polylogarithmic currents can be obtained by gluing the objects resulting to the local construction we are going to explain.\\

\subsection*{A1.1 The fibrewise exponential map} 

We have the following exact sequence of abelian groups over $S$:
\begin{eqnarray} \label{exactseq}
0 \to \Lambda \overset{i}{\to} e^* \text{T} (X^{\infty} / S^{\infty} ) \overset{\exp_{X^{\infty}  / S^{\infty} }}{\to} X^{\infty} \to 0,
\end{eqnarray}
where $\text{T} (X^{\infty} / S^{\infty} )$ is the relative tangent bundle of $\pi \colon X^{\infty} \to S^{\infty}$ and 
$\exp_{X^{\infty}  / S^{\infty}}$ is the fibrewise exponential map. The monomorphism $i$ induces an isomorphism
$ \widetilde{i} \colon \Lambda \otimes \OO_{S^{\infty}} \isom e^* \text{T} (X^{\infty} / S^{\infty} )$.\\

\subsection*{A1.2 $\mathcal{C}^{\infty}$-trivialisation} 

\noindent {\bf Assumption} $-$ Let $S$ be  simply connected. \\

We fix a base point $s_0 \in S$.
Since $S$ is simply connected, there exists a canonical isomorphism $ \iota \colon \underline{\Lambda_{s_0}} \isom \underline{\Lambda}$,
where $\underline{\Lambda_{s_0}}$ is the constant sheaf on $S$ associated to $\Lambda_{s_0} = H_1(A_{s_0} , \Z )$.
The isomorphism $\widetilde{i} \circ \iota \otimes {\text Id}_{ \OO_{S^{\infty}} }    \colon
\underline{\Lambda_{s_0}} \otimes \OO_{S^{\infty}} \to e^* \text{T} (X^{\infty} / S^{\infty} )$
of sheaves of locally free $\OO_{S^{\infty}}$-modules 
corresponds to a unique isomorphism of vector bundles over $S^{\infty}$ which we also denote 
$\widetilde{i} \circ \iota \otimes {\text Id}_{ \OO_{S^{\infty}} }$
$$ \widetilde{i} \circ \iota \otimes {\text Id}_{ \OO_{S^{\infty}} } \colon S^{\infty} \times 
(\Lambda_{s_0} \otimes \R) \to e^* \text{T} (X^{\infty} / S^{\infty} ).$$
We observe  that 
$\exp_{X^{\infty}  / S^{\infty} } \circ (\widetilde{i} \circ \iota \otimes {\text Id}_{ \OO_{S^{\infty}} }) \colon 
S^{\infty} \times 
(\Lambda_{s_0} \otimes \R) \to X^{\infty}$ induces (cf exact sequence (\ref{exactseq})) an isomorphism
of families of real tori over $S^{\infty}$ 
$$  \varphi \colon 
S^{\infty} \times  (\Lambda_{s_0} \otimes \R   /\Lambda_{s_0})   \to X^{\infty}.$$

The tangent bundle of $\Lambda_{s_0} \otimes \R   /\Lambda_{s_0}$ can be naturally identified with the trivial vector bundle
$(\Lambda_{s_0} \otimes \R   /\Lambda_{s_0}) \times  ( \Lambda_{s_0} \otimes \R)$.
Thus we have a natural identification
\begin{eqnarray} \label{tangentbdl}
\text{T} (S^{\infty} \times  (\Lambda_{s_0} \otimes \R   /\Lambda_{s_0})) =  
 pr_1 ^* \;  \text{T} S^{\infty} \oplus (S^{\infty} \times  (\Lambda_{s_0} \otimes \R   /\Lambda_{s_0})) \times  ( \Lambda_{s_0} \otimes \R),
\end{eqnarray}
 where $pr_1$ is the canonical projection  $S^{\infty} \times  (\Lambda_{s_0} \otimes \R   /\Lambda_{s_0}) \to S^{\infty}.$

\subsection*{A1.3 Polarisation form and symplectic pairings}

The $(1,1)$-form $\omega$ on $X^{\infty}$ induces a pairing $< \cdot , \cdot > \colon \Lambda_{s_0} \wedge \Lambda_{s_0} \to \Z(1)$.
Extending this pairing by linearity we get two other pairings
$$ (\Lambda_{s_0} \otimes \C) \wedge (\Lambda_{s_0} \otimes \C) \to \C, \quad 
\quad \text{ and } \quad (\Lambda_{s_{0}} \otimes \OO_{X^{\infty}}) \wedge (\Lambda_{s_{0}} \otimes \OO_{X^{\infty}}) \to \OO_{X^{\infty}}$$
also denoted by the symbol $< \cdot , \cdot >$. \\

We remark that 
the symplectic pairing  $< \cdot , \cdot > \colon (\Lambda_{s_0} \otimes \C) \wedge (\Lambda_{s_0} \otimes \C) \to \C$, 
which corresponds to an element in $\bigwedge^2 \Hom_{\C}(\Lambda_{s_0} \otimes \C , \C)$, 
induces a complex differentiable 2-form on $S^{\infty} \times  (\Lambda_{s_0} \otimes \R   /\Lambda_{s_0})$ (see (\ref{tangentbdl}))
which corresponds to $2 \; \varphi^* \omega$ (cf proof of the proposition 2.2.4 in \cite{l}).

\subsection*{A1.4 Definition of the functions $\chi_{\lambda}$ }

Let $\Lambda'_{s_0}$ be the $2 \pi i$-dual of $\Lambda_{s_0}$ with respect to $< \cdot , \cdot >$, i.e.
$$\Lambda'_{s_0} := \{ \lambda' \in \Lambda_{s_0} \otimes \C \; | \; < \lambda' , \lambda > \in 2 \pi i \Z \}, $$
$\kappa$ be the index $[ \Lambda'_{s_0} : \Lambda_{s_0}]$ and  $\lambda \in \Lambda'_{s_0}$.
 We define a complex valued function $\chi_{\lambda}$ on $X^{\infty}$ by
$$ \chi_{\lambda} ( x ) = \exp ( < \lambda ,    pr_2 \circ \varphi^{-1} ( x )   > ) \quad \text{ for all }  x \in X,$$
where $pr_2$ denotes the natural projection 
 $S^{\infty} \times ( \Lambda_{s_0} \otimes \R   /\Lambda_{s_0}) \to \Lambda_{s_0} \otimes \R   /\Lambda_{s_0}.$

\subsection*{A1.5 Construction of vector fields}

The tangent map of $\varphi$ induces an isomorphism (see (\ref{tangentbdl}))
$$ pr_1 ^* \;  \text{T} S^{\infty} \otimes \C \oplus (S^{\infty} \times 
 (\Lambda_{s_0} \otimes \R   /\Lambda_{s_0})) \times  ( \Lambda_{s_0} \otimes \C) 
\isom \text{T} A^{\infty} \otimes \C. $$

Using this isomorphism we associate to each section  of 
\begin{eqnarray} \label{hh}
 \hh :=  \Lambda \otimes \C \otimes \OO_{S^{\infty}}
 \overset{ \iota^{-1} \otimes \text{Id}_{\C \otimes \OO_{S^{\infty  } } } }{  =}
\Lambda_{s_0}  \otimes \C  \otimes \OO_{S^{\infty}}
\end{eqnarray}
a complex vertical vector field on $X^{\infty}$.\\

\noindent {\bf Convention} $-$
We can canonically associate to a section of $\hh$ over $S^{\infty}$ a section of $\pi^* \hh$ over $X^{\infty}$
and 
a  complex vertical vector field on $X^{\infty}$ as explained before. 
These three objects are denoted by the same symbol. 
For example, if $\lambda \in \Lambda'_{s_0}$ (viewed as a section of $\hh$ over $S^{\infty}$ via (\ref{hh})) 
and $\lambda = \lambda^{-1,0} + \lambda^{0,-1}$
is the decomposition of $\lambda$ with respect to the Hodge decomposition
$ \hh = \hh^{-1,0} \oplus \hh^{0,-1}$, the convention holds for the sections $\lambda^{-1,0}$ and $\lambda^{0,-1}$
of $\hh$.

\subsection*{A1.6 Definition of the polylogarithmic current}

The polylogarithmic current is defined as 
$g := \displaystyle  \sum_{n \geq 2} \quad \sum_{a,b \geq 1, \; a+b=n}   (-1)^a g_{a,b}$
where for each $a,b\geq~1$ 

$$
g_{a,b} := \displaystyle \frac{(-1)^d}{d! \kappa}  \; 
\sum_{\lambda \in \Lambda_{s_0}' \setminus \{ 0\} } 
i_{\lambda^{-1,0}} i_{\lambda^{0,-1}} \chi _{\lambda}
 \frac{(\lambda^{0,-1})^{a-1}
(\lambda^{-1,0})^{b-1}}{ (<\lambda^{-1,0},\lambda^{0,-1}> -
[\lambda^{-1,0}] )^{a+b}} \; \omega^d.
$$
Here $i_{\lambda^{-1,0}}$ (resp. $i_{\lambda^{0,-1}}$) denotes the contraction operator associated to the vector field 
$\lambda^{-1,0}$ (resp. $\lambda^{0,-1}$) and $[\lambda^{-1,0}]$ is the Lie derative corresponding to the 
vector field  $\lambda^{-1,0}$.
Using the power series expansion of $(c-x)^{-(a + b)}$  and the vanishing of
$[\gamma^{-1,0}]^k \omega^d$ for $k > 2d$ \cite[Prop 3.2.2]{l}, we get

$$ g_{a,b} = 
\sum_{k=0}^{2d} \frac{ (-1)^d (a+b +k-1)!}{(a+b-1)!k!d! \kappa}
\underset{=: g_{a,b}^k }{\underbrace{\sum_{\lambda \in  \Lambda_{b_0}' \setminus \{ 0\} }
 \chi _{\lambda}
 \frac{(\lambda^{0,-1})^{a-1}
(\lambda^{-1,0})^{b-1}} { (<\lambda^{-1,0},\lambda^{0,-1}>
)^{a+b+k}}i_{\lambda^{-1,0}}i_{\lambda^{0,-1}}
[\lambda^{-1,0}]^k \; \omega^d }}.$$

One can check that the definition of $g$ does not depend on the choice of the base point $s_0$.

\section*{A2. A smoothness result for the polylogarithmic current}

We keep the notations of the previous part. \\

\noindent 
{\bf Proposition A2.1} $-$ {\it The restriction of the polylogarithmic current $g$ over $X^{\infty} \setminus e(S^{\infty})$ is a smooth current.} \\

\noindent {\it Proof} $-$ 
Since the smoothness is a local property we may assume that $S$ is simply connected and there
exists global coordinates $x_1, \dots, x_{2r}$ on $S^{\infty}$, where $r$ is the dimension of the complex analytic variety $S$. 
We fix a base point $s_0 \in S$ as in the part A1.
Considering the definition of $g$ recalled in the part A1 we observe that it is enough to prove the smoothness of $g_{a,b}^k$  
over $X^{\infty} \setminus e(S^{\infty})$  ($a, b \geq 1$, $0 \leq k \leq 2d$). This is equivalent to prove the smoothness 
of $\varphi^*g_{a,b}^k$ on $S^{\infty} \times ( \Lambda_{s_0} \otimes \R   / \Lambda_{s_0} \setminus \{0 \} )$. \\

We observe that $\varphi^* g_{a,b}^k$ is an expression of the shape
$$ \sum_{\lambda \in \Lambda_{s_0}' \setminus \{ 0 \} } \chi' _{\lambda}
\frac{P(\lambda)} { (Q(\lambda) )^{m}}$$ 
where 
$\chi'_{\lambda} \colon S^{\infty} \times (\Lambda_{s_0} \otimes \R   /\Lambda_{s_0}) \to \C$ is the smooth function 
defined by $\chi'_{\lambda}(s,u) = \exp( < \lambda , u > )$ for all $(s,u) \in S^{\infty} \times (\Lambda_{s_0} \otimes \R   /\Lambda_{s_0})$,
 $Q$ is a positively definite quadratic form on the lattice $\Lambda_{s_0}$ and $P$ is a
homogeneous polynomial function of degree $m$ on the lattice with
values in some finite-dimensional vector space (this space is
tensor product of three spaces: 
the symmetric power  $\Sym^k (\Lambda_{s_0} \otimes \C)$, 
the space $\bigwedge^\bullet  \Hom_{\C}(\Lambda_{s_0} \otimes \C , \C)$ 
and the space $\bigwedge^\bullet \text{Span}_{\C}(dx_1, \dots , dx_{2r})$.
The quadratic form and the polynomial are smooth
functions on the base $S^{\infty}$.\\

Consider the convergent for $s \in \C$ such that $\Re(s)\gg 0$ series ${\cal K}(Q,P,s)$
$$ \sum_{\lambda \in \Lambda_{s_0}' \setminus \{ 0 \} }  \chi' _{\lambda}
\frac{P(\lambda)} { (Q(\lambda) )^{s}}.$$
It converges uniformly with respect to the base. For rather big $\Re(s)$ we have the following expression for the
product $\Gamma(s){\cal K}(Q,P,s)$ of ${\cal K}(Q,P,s)$ and the
$\Gamma$-function ($\displaystyle \Gamma(s)=\int_0^\infty t^{s-1}e^{-t}d\,t$) :

$$\Gamma(s){\cal K}(Q,P,s) =
\sum_{ \lambda \in \Lambda_{s_0}' \setminus \{ 0 \}  } \chi'_{\lambda}
\Gamma(s)\frac{P(\lambda)} { (Q(\lambda) )^{s}}= \int_0^\infty
\sum_{ \lambda \in \Lambda_{s_0}' \setminus \{ 0 \} } \chi'_{\lambda}
 t^{s-1}e^{-t}\frac{P(\lambda)} { (Q(\lambda) )^{s}}d\,t
$$
$$=\int_0^\infty
\sum_{\lambda \in \Lambda_{s_0}' \setminus \{ 0 \}} \chi'_{\lambda}
 e^{-t(Q(\lambda) )} {P(\lambda)} t^{s-1} d\,t=
\int_0^\infty (\Theta(Q,P,u,t)  -P(0))t^{s-1} d\,t,
$$
where 
$\Theta(Q,P,u,t)= \displaystyle
\sum_{\lambda \in \Lambda_{s_0}' } \exp(< \lambda , u > )  e^{-t(Q(\lambda) )} {P(\lambda)}$ for
$u \in \Lambda_{s_0} \otimes \R   /\Lambda_{s_0}$, $t \in \R ^{>0}$. \\

We split the domain of integration into  two subdomains: from zero
to some nonzero constant $A$ and from $A$ to infinity.\\

Over a compact subset on the base we have a bound 
$Q(\lambda) > C \sum |\lambda_i|^2$ where $\lambda_i$ are coordinates of $\lambda$ with
respect to some basis of $\Lambda_{s_0}'$ and $C$ is  some positive real number.
Hence for $t\gg 0$ one have an uniform with respect to the base $S$ bound $(\Theta(Q,P,u,t)
 -P(0))=O(\exp(-Kt))$, so the integral from
 $A$ to $\infty$ converges for any $s$.\\
 
The integral from $0$ to $A$ can be calculated via Poisson summation formula. 
We recall this formula in our context.
Let $\text{vol}$ be a volume form on $\Lambda_{s_0} \otimes \R$ such that the covolume
of $\Lambda_{s_0}$ equals 1. 
For a rapidly decreasing
function $f$ on  $\Lambda_{s_0}$ denote by $\tilde{f}$ its Fourier
transform with respect to the pairing 
$< \cdot , \cdot  > \colon (\Lambda_{s_0} \otimes \R) \wedge (\Lambda_{s_0} \otimes \R) \to \R(1)$
 and the volume form $\text{vol}$:
$$\tilde{f}(p)=\int_{\Lambda_{s_0} \otimes \R   }f(x)\exp(<x,p>) \text{vol}_x, \quad p \in \Lambda_{s_0} \otimes \R .$$
Then  $\displaystyle \sum_{\lambda' \in \Lambda_{s_0}' }f(\lambda')=\sum_{\lambda \in \Lambda_{s_0}} \tilde{f}(\lambda)$. \\
\\

Let $h \in  \Lambda_{s_0} \otimes \R$ and 
let $u$ be the image of $h$ under the natural projection $\Lambda_{s_0} \otimes \R \to \Lambda_{s_0} \otimes \R   /\Lambda_{s_0}.$\\

The value at $p \in \Lambda_{s_0} \otimes \R$ of the Fourier transform of the function 
$ \exp (< x, h >) \exp( - t Q(x) ) {P(x)}$ in $x \in \Lambda_{s_0} \otimes \R$ is equal to
$$
\pi^d {\rm Disc}(tQ)^{-1/2}\exp\left(-\frac{\pi^2}{t}Q^{\vee}(p+h)\right)\hat{P}(p)
 =
t^{-d}\pi^d {\rm Disc}(Q)^{-1/2}\exp\left(-\frac{\pi^2}{t}Q^{\vee}(p+h)\right)\hat{P}(p),
$$
where ${\rm Disc}$ denotes the discriminant of the quadratic form with respect to the volume form $\text{vol}$, 
$Q^{\vee}$ is the dual (with respect to the  pairing $< \cdot , \cdot >$) to $Q$ quadratic form and
$\hat{P} $ is some polynomial of the same degree as $P$.\\

We denote the sum $ \displaystyle \sum_{\lambda \in \Lambda_{s_0}} e^{-t (\pi^2 Q^{\vee}(\lambda+h ))} {\hat{P}(\lambda)}$
(which depends only on $u$) by $\hat{\Theta}(Q^{\vee},\hat{P},u,t)$.\\

 Then from  the Poisson summation formula we get
$$\Theta(Q,P,u,t)= t^{-d}\pi^d {\rm Disc}(Q)^{-1/2}\hat{\Theta}(Q^{\vee},\hat{P},u,t^{-1}).$$ 
Hence 
$$
\begin{array}{lll} 
\displaystyle \int_0^A(\Theta(Q,P,u,t)  -P(0)) t^{s-1} d\,t & =  &
\displaystyle \pi^d {\rm Disc}(Q)^{-1/2} \int_0^At^{-d} \hat{\Theta} (Q^{\vee},\hat{P},u,t^{-1}) t^{s-1} d\,t -P(0) \int_0^A t^{s-1} d\,t  \\
&  & \\
& = &
 \displaystyle \pi^d {\rm Disc} (Q)^{-1/2} \int_{A^{-1}}^{\infty} \hat{\Theta} (Q^{\vee},\hat{P},u,x) x^{d-s-1} d\,x -P(0) \frac{1}{s} t^s |_0^A .\\
\end{array}
$$

For $x\gg 0$ one has an uniform with respect to the base $S$ and $u \neq 0$ bound 
$\hat{\Theta}(Q^{\vee},\hat{P},u,x)  =O(\exp(-Kx)).$
So the first summand convergent integral is a smooth function in
$u\neq 0$ for any $s$. The second $-P(0)A^s/s$ vanishes as $P$ is
homogeneous. This finishes the proof of the smoothness of the
polylogarithmic currents.\\

\hfill{$\Box$}


\newpage
\begin{thebibliography}{99}

\bibitem[BBD]{bbd} A.A. Beilinson, J. Bernstein, P. Deligne, \emph{Faisceaux pervers},
Analyse et topologie sur les espaces singuliers, Astérique 100.




\bibitem[BL]{bl}
  A.A. Beilinson, A. Levin, \emph{The Elliptic Polylogarithm}, in U. Jannsen,
  S.L. Kleiman, J.P. Serre, ``Motives'', Proceedings of the research Conference on Motives
  held July 20 - August 2, 1991, in Seattle, Washington,
  Proc. of Symp. in Pure Math. 55, Part II, AMS, p. 123--190 (1994).  


\bibitem[B]{b} D. Blotti\`ere, \emph{Les classes d'Eisenstein des vari\'et\'es de Hilbert-Blumenthal},
               en pr\'eparation.











\bibitem[D]{d} J. Dieudonné, \emph{\'Eléments d'analyse}, Tome III, \'Editions Jacques Gabay, 2003.





\bibitem[GH]{gh} P. Griffiths, J. Harris, \emph{Principles of algebraic geometry}, Wiley Classics Library.


%






\bibitem[HZ]{hz}
  R. M. Hain, S. Zucker,
  \emph{Unipotent variations in mixed Hodge structures},
 Inv. Math. 88, p. 83--124 (1987).

\bibitem[I]{i} B. Iversen, {\it Cohomology of sheaves}, Universitext, Springer-Verlag.

\bibitem[Ka]{ka}
  M. Kashiwara, \emph{A study of variations of mixed Hodge structure}, Publ. RIMS,
  Kyoto Univ. 22,
  p. 991--1024 (1986).




\bibitem[Ki]{k}
  G. Kings, \emph{$K$-theory elements for the polylogarithm of abelian schemes},
  J. reine angew. Math. 517, p. 103--116 (1999).





\bibitem[L]{l}
  A. Levin, \emph{Polylogarithmic Currents on Abelian Varieties}, 
  in A. Reznikov, N. Schappacher (eds),
  ``Regulators in Analysis, Geometry and Number Theory'', 
  Proc. of the German-Israeli workshop held
  March 11--20 1996 at Landau Center, Jerusalem, 
  Progress in Math. 171, Birkhäuser, p. 207--229 (2000).



   







\bibitem[S]{s}
  M. Saito, \emph{Mixed Hodge Modules}, Publ. RIMS Kyoto Univ. 26, p. 221--233 (2000).

  





\bibitem[T]{t}
A. Terras, 
\emph{Harmonic analysis on symmetric spaces and applications I},
 Springer-Verlag New York (1985).

\bibitem[W]{w}
  J. Wildeshaus, \emph{Realizations of polylogarithms}, L.N.M. 1650, Springer-Verlag Berlin (1997).




\end{thebibliography}
\end{document}